\documentclass{tran-l}

\usepackage{amssymb}

\newtheorem{theorem}{Theorem}[section]
\newtheorem{lemma}[theorem]{Lemma}
\newtheorem{corollary}[theorem]{Corollary}
\newtheorem{proposition}[theorem]{Proposition}

\theoremstyle{definition}
\newtheorem{definition}[theorem]{Definition}

\theoremstyle{remark}

\numberwithin{equation}{section}

\newcommand\Aut[1]{\text{Aut}(#1)}

\def\bec{\begin{corollary}}
\def\enc{\end{corollary}}

\def\bet{\begin{thm}}
\def\ent{\end{thm}}

\def\becor{\begin{corollary}}
\def\encor{\end{corollary}}

\def\bel{\begin{lem}}
\def\enl{\end{lem}}

\def\bedef{\begin{definition}}
\def\endef{\end{definition}}

\def\bep{\begin{proposition}}
\def\enp{\end{proposition}}

\newcommand{\comment}[1]{}

\newcommand{\Fld}{\textbf{Fld}}

\newcommand{\fin}{\not\cong}
\newcommand{\Inf}{\textbf{Inf}}

\newcommand{\set}[2]{\ensuremath{ \{ #1 : #2 \} }}

\renewcommand{\deg}[1]{\ensuremath{\text{deg}(#1)}}

\newcommand{\Z}{\mathbb{Z}}
\newcommand{\F}{\mathbb{F}}
\newcommand{\Q}{\mathbb{Q}}
\newcommand{\R}{\mathcal{R}}
\newcommand{\C}{\mathcal{C}}
\newcommand{\A}{\mathcal{A}}
\newcommand{\B}{\mathcal{B}}

\newcommand{\fbar}{\overline{f}}
\newcommand{\gbar}{\overline{g}}

\newcommand{\Qbar}{\overline{\mathbb{Q}}}

\newcommand{\Ftilde}{\widetilde{F}}

\newcommand{\Ltilde}{\widetilde{L}}

\newcommand{\ptilde}{\tilde{p}}
\newcommand{\qtilde}{\tilde{q}}

\newcommand{\utilde}{\tilde{u}}
\newcommand{\vtilde}{\tilde{v}}

\newcommand{\xtilde}{\tilde{x}}
\newcommand{\ytilde}{\tilde{y}}
\newcommand{\ztilde}{\tilde{z}}

\newcommand{\Qtilde}{\tilde{Q}}

\renewcommand{\L}{\mathcal{L}}
\newcommand{\M}{\mathfrak{M}}

\newcommand{\avec}{\vec{a}}
\newcommand{\bvec}{\vec{b}}

\newcommand{\vvec}{\vec{v}}

\newcommand{\xvec}{\vec{x}}
\newcommand{\Xvec}{\vec{X}}

\newcommand{\cS}{{\mathcal S}}

\newcommand{\Gal}[2]{\text{Gal}(#1/#2)}
\newcommand{\la}{\langle}
\newcommand{\ra}{\rangle}

\newcommand{\sqrtp}{\sqrt{p}}

\def\converges{\!\downarrow}

\def\Ss{{\mathfrak S}}

\newcommand{\rg}[1]{\text{range}(#1)}

\newcommand{\bfd}{\boldsymbol{d}}
\def\bfz{\boldsymbol{0}}
\def\s01{\ensuremath{\Sigma^0_1}}
\def\d02{\ensuremath{\Delta^0_2}}
\def\phi{\varphi}

\def\res{\!\!\upharpoonright\!}

\begin{document}

\title[Categoricity Properties for Computable Algebraic Fields]{Categoricity Properties for\\Computable Algebraic Fields}


\author[D.R.\ Hirschfeldt]{Denis R.\ Hirschfeldt}
\address{Department of Mathematics \\
University of Chicago\\
5734 S.\ University Ave.\\
Chicago, IL 60637 U.S.A.}
\email{drh@math.uchicago.edu}
\thanks{The first author was partially supported by Grants \# DMS--0801033 and
DMS--1101458 from the National Science Foundation.
The second author was partially supported by Grant \# DMS--0739346 from the National Science Foundation.
The third author was partially supported by Grant \# DMS--1001306 from
the National Science Foundation, by Grant \# 13397 from the
Templeton Foundation, by the Centre de Recerca Matem\'atica, 
and by several grants from The City University of New York 
PSC-CUNY Research Award Program.
The fourth author was partially supported by Grants \# DMS--0650927
and DMS--1161456
from the National Science Foundation, by Grant \# 13419 from the
Templeton Foundation, and by an ECU Faculty Senate Summer 2011 Grant.}

\author[K.\ Kramer]{Ken Kramer}
\address{Department of Mathematics \\
Queens College -- C.U.N.Y. \\
65-30 Kissena Blvd.\\
Flushing, New York  11367  U.S.A.;
Ph.D. Program in Mathematics \\
C.U.N.Y.\ Graduate Center \\
365 Fifth Avenue \\
New York, NY  10016 U.S.A.}
\email{kkramer@qc.cuny.edu}

\author[R.\ Miller]{Russell Miller}
\address{Department of Mathematics \\
Queens College -- C.U.N.Y. \\
65-30 Kissena Blvd.\\
Flushing, New York  11367  U.S.A.;
Ph.D. Programs in Mathematics  and Computer Science\\
C.U.N.Y.\ Graduate Center \\
365 Fifth Avenue \\
New York, NY  10016 U.S.A.}
\email{Russell.Miller@qc.cuny.edu}

\author[A.\ Shlapentokh]{Alexandra Shlapentokh}
\address{East Carolina University\\
Department of Mathematics \\
Greenville, NC 27858 U.S.A.}
\email{shlapentokha@ecu.edu}




\begin{abstract}
We examine categoricity issues for computable algebraic fields.
We give a structural criterion for relative computable
categoricity of these fields, and use it to construct a field that is
computably categorical, but not relatively computably categorical.
Finally, we show that computable categoricity for this class of fields
is $\Pi^0_4$-complete.
\end{abstract}

\maketitle

\section{Introduction}
\label{sec:intro}

Fields were the first class of structures for which the notion
of computable categoricity was ever expressed.  In their landmark
study of effectiveness in field theory, Fr\"ohlich and
Shepherdson presented ``two explicit fields which are
isomorphic but not explicitly isomorphic'' \cite[Corollary 5.51]{FS56}.
In modern terminology, we would say that these two fields
are both computable, and are classically isomorphic
but not computably isomorphic.  Thus they fail to satisfy
the definition of computable categoricity.
\begin{definition}
\label{defn:cc}
The \emph{Turing degree} of a countable structure $\A$ is the join of the 
degrees
of the functions and relations of $\A$, or equivalently, the Turing
degree of its atomic diagram. A \emph{computable structure} is one with 
Turing degree $\bf 0$.

A computable structure $\A$ is \emph{computably categorical}
if, for every computable structure $\B$ isomorphic
to $\A$, there exists a computable isomorphism from $\A$ onto $\B$.

More generally, a computable structure $\A$ is \emph{relatively
computably categorical}
if, for every structure $\B$ with domain $\omega$ that is
isomorphic to $\A$, there exists an isomorphism from $\A$ onto $\B$
that is computable in the Turing degree of $\B$.
\end{definition}
For these and other definitions from
computable model theory, \cite{AK00} and \cite{H98} are excellent sources.
The article \cite{MS11}, written by two of the present authors, also serves
to introduce these and many related concepts in more detail, and the articles
\cite{Notices08} and \cite{Proceedings08} present the basic notions about
computable fields for readers unfamiliar with them.

For over fifty years since that first result in \cite{FS56},
computable categoricity for fields has remained largely
a mystery.  For many other classes of structures, mathematicians
have found structural definitions equivalent to computable categoricity:
see for instance \cite{GD80}, \cite{GLS03}, \cite{KKM05},
\cite{LMMS05}, \cite{M05}, \cite{R81a}, and \cite{R81b}.
As an example, Goncharov and Dzgoev, and independently Remmel,
showed that a linear order is computably categorical if and only if it has
only finitely many pairs of adjacent elements. (Two distinct elements
of a linear order are \emph{adjacent} if there is no element of the
order between them.) This criterion
is not quite expressible in first-order model theory, since
it involves finiteness, but intuitively it is distinctly more
``structural'' than Definition \ref{defn:cc}.  In terms of computational
complexity this criterion is $\Sigma^0_3$ (and is readily shown
to be complete at that level), whereas the statement
of Definition \ref{defn:cc} is $\Pi^1_1$, quantifying
over all possible (classical) isomorphisms.  Indeed, for linear orders,
computable categoricity turns out to coincide with relative
computable categoricity, and Ash, Knight, Manasse, and Slaman
established in \cite{AKMS89} that relative computable categoricity
is always a $\Sigma^0_3$ property.  (Unpublished work \cite{Chis}
by Chisholm yields the same result.)  On the other hand,
although relative computable categoricity clearly
implies computable categoricity, it was established
independently by Khoussainov and Shore in \cite{KS97}
and by Kudinov in \cite{K96} that the two concepts are
not equivalent.  More recently, Downey,
Hirschfeldt, and Khoussainov showed in \cite{DHK03} that
relative computable categoricity can be viewed as a kind of uniform
version of computable categoricity, although this fact was already 
implicit in work of Ventsov \cite{V92}.

For fields, however, only a few significant criteria for computable
categoricity (or for its failure) have been discovered.  The situation
is straightforward when the field is algebraically closed:  Ershov showed in \cite{E77}
that such a field is computably categorical if and only if it has finite transcendence degree
over its prime subfield (either $\Q$ or the $p$-element field $\F_p$,
depending on characteristic).  Earlier, Fr\"ohlich and Shepherdson \cite{FS56}
had established that all normal algebraic extensions of $\Q$ and
of $\F_p$ are computably categorical.  These results failed to extend
to fields more generally, however:  algebraic extensions of $\Q$
that are not computably categorical have been known at least since 
\cite{E77},
and Miller and Schoutens recently constructed a computably categorical
field of infinite transcendence degree over $\Q$ (see \cite{MS08}).

The transcendence degree of the field over its prime subfield
is soon seen to be of paramount importance in these considerations.
For algebraic field extensions $F$ of $\Q$, one can identify each element
$x\in F$ to within finitely many possibilities by finding the minimal polynomial
of $x$ in $\Q[X]$, and likewise for algebraic extensions of $\F_p$;
this fact follows from the existence of \emph{splitting algorithms} for 
$\Q$ and for each $\F_p$.  When one wishes to compute an isomorphism
between two such fields, the task of determining an image for $x$
is not completely solved by this knowledge, but its degree of difficulty
becomes relatively low; see \cite{M09b} for the current state
of knowledge on this topic.  The paper \cite{MS11},
written by two of us, is in many ways a precursor to this paper,
and produces a criterion for computable categoricity in case
the entire algebraic field $F$ has a splitting algorithm:  such an $F$
is computably categorical if and only if its orbit relation is decidable,
in which case it is also relatively computably categorical.
(The \emph{orbit relation} holds of the pair $\la a,b\ra\in F^2$
if and only if some automorphism of $F$ maps $a$ to $b$.)  Below
we prove that this criterion does not extend to all computable algebraic
fields; indeed both implications fail.  (For further background about
splitting algorithms and related concepts, including the
\emph{splitting set} and the \emph{root set} of $F$,
we suggest \cite{Notices08}, \cite{M09a}, and \cite{MS11}.)

This paper focuses on computable algebraic fields $F$.  We do not
assume the existence of a splitting algorithm for $F$, although
our results do apply in the situation where $F$ has a splitting algorithm.
That case was mostly explained in \cite{MS11}, however,
while here we show that the situation without a splitting algorithm
is significantly more difficult.  In particular, computable algebraic
fields without splitting algorithms can be computably categorical
without being relatively computable categorical (see Theorem \ref{thm:notrcc}).
Moreover, the complexity of computable categoricity goes up
when the field is not required to have a splitting algorithm:
computable categoricity is $\Pi^0_4$-complete for algebraic fields
(see Theorem \ref{thm:Pi04complete}),
whereas with a splitting algorithm it is equivalent to relative computable
categoricity, hence only $\Sigma^0_3$-complete.  The increase
in complexity is significant, but the switch from $\Sigma$ to $\Pi$
is also significant.  Indeed, for algebraic fields, Definition \ref{defn:cc}
has complexity $\Pi^0_4$, since the property of being isomorphic
is only $\Pi^0_2$, distinctly simpler than the usual $\Sigma^1_1$.
(This fact is based on Corollary \ref{cor:nofiniteblock}.)
Therefore, our results show that the standard definition of
computable categoricity actually has the minimum possible 
complexity all by itself, when restricted to algebraic fields:
no structural (or other) criterion can improve it.  To our knowledge,
algebraic fields are the first class of structures for which
this has been shown to be the case.

\comment{
We do show algebraic fields to be nice structures in one related respect.
Goncharov defined the \emph{computable dimension} of a computable
structure to be the number of computable presentations of that structure,
up to computable isomorphism.  He showed that every cardinal
from $1$ through $\omega$ can be the computable dimension of a computable
structure.  (See \cite{G82} and \cite{G98} for these and related results.)
However, by far the most common computable dimensions
are $1$ (which is equivalent to computable categoricity) and $\omega$,
and for many classes of structures, these are the only
possible computable dimensions:  linear orders, Boolean algebras,
and trees, for example.  We show in Corollary \ref{cor:dimomega}
that algebraic fields too can only have computable dimension
$1$ and $\omega$.  Since the intermediate dimensions are
usually regarded as pathological, this fact makes algebraic fields
seem like a nice class of structures, whereas the result about
computably categorical fields failing to be relatively computably categorical
suggested the opposite.
Indeed, we believe that algebraic fields constitute the first known example
of a class of structures (or at least, a class commonly seen in mathematics)
for which $1$ and $\omega$ are the only possible computable dimensions,
yet computable categoricity does not imply relative computable categoricity.
So our results differentiate these two pathologies from each other.
}

\section{Useful Results on Computable Fields}
\label{sec:fields}

Substantial work on computable algebraic fields and categoricity
has appeared recently, giving rise to several useful techniques
for constructing computable fields.  In this section we review
assorted properties of algebraic fields relevant to these techniques,
with references to allow the reader to look up their proofs and
to see how they were originally used.

The following result, which appears as Lemma 2.10 in \cite{M09b},
will often save us from having to worry about surjectivity
as we compute isomorphisms between fields.
\begin{lemma}
\label{lemma:endo}
For an algebraic field $F$, every endomorphism (i.e.\ every
injective homomorphism $g:F\to F$) is an automorphism.
\end{lemma}
\begin{corollary}
\label{cor:endo}
If $E\cong F$ are isomorphic algebraic fields, and $f:E\to F$
is a field embedding (by which we mean a field homomorphism
with $f(1)\neq 0$), then the image of $f$ is all of $F$.
That is, such an $f$ must be an isomorphism.
\end{corollary}

We will use the standard notation for Galois groups:
if $F\subseteq K$ is a \emph{Galois extension}
(i.e.\ an algebraic normal separable field extension), then $\Gal{K}{F}$,
the \emph{Galois group} of $K$ over $F$, is the
group of all field automorphisms of $K$ which
restrict to the identity map on $F$.
As we build computable fields, it frequently happens that,
having already built a computable field $F_s$,
we wait to see whether a particular function
will converge on a particular input.  If it does not converge,
then $F_s$ itself satisfies a particular requirement $\R_2$ for
the construction, whereas if it does converge, we can add
more elements to $F_s$ to build the larger field $K_2$ and
satisfy the requirement that way.  When considering two
distinct requirements, it is useful to be certain
that extending $F_s$ to $K_2$ to satisfy $\R_2$ will not
disrupt our plan to build a different extension $K_1$
if necessary to satisfy a different requirement $\R_1$.
Usually, if $\Gal{K_1}{F_s}\cong\Gal{E}{K_2}$
(where $E$ is the field generated by $K_1$ and $K_2$ together),
we can avoid the disruption to $\R_1$, and one way
to ensure this isomorphism holds is to make $K_1\cap K_2=F_s$.
(See \cite[p.\ 243, Exercise 2]{J85}, for example.)

To achieve this end, we will often use a Galois extension
of $\Q$ whose Galois group over $\Q$ is
the symmetric group on the roots of a given polynomial, since
this choice allows us to adjoin some of these roots immediately and
keep others out of the field until needed.  (The proof of Theorem
\ref{thm:rccorbit} is a good example of such a construction.)
Therefore, it is frequently useful for an extension
such as the $K_i$ above to have symmetric Galois group
over the current ground field, as this property ensures that it is the
splitting field of a polynomial whose roots are essentially
all independent of each other.  
The first theorem for this purpose appeared as Theorem 2.15 in \cite{M09a},
and provides a supply of such extensions.  The proof given there
was devised by Kevin Keating.  Since we also want these extensions
not to interfere with each other (and since extensions with large
symmetric Galois groups cannot be taken to have relatively prime degrees,
which is the most obvious way to avoid such interference),
we now extend that theorem to include the linear disjointness
of the extensions.

\begin{definition}
\label{defn:lindisjt}
Two Galois extensions $E\subseteq K$ and $E\subseteq L$ within a larger
field $F$ are \emph{linearly disjoint} if $K\cap L=E$.
(This is a particular case of the definition of linear disjointness
for algebraic field extensions in general, which requires that
$K$ and $L$ together generate an extension whose degree
over $E$ is the product $[K:E]\cdot [L:E]$.)
\end{definition}
This means that we can add elements of $K$ to $E$ to build $F$
to satisfy one requirement, and close under the field operations, without worrying
that these new elements might accidentally force certain elements of $L$
to enter $F$ as well and thereby upset our satisfaction of a different requirement.
The simplest case of linear disjointness occurs when the degrees
$[K:E]$ and $[L:E]$ are relatively prime:  the degree $[K\cap L:E]$
divides both, hence equals $1$, so $K\cap L= E$.  (Indeed, in this
situation $K$ and $L$ need not be Galois extensions of $E$.)

Proposition \ref{prop:Ken} gives the recursive step for our procedure
for building many distinct extensions of $\Q$, each one linearly
disjoint from the field generated by all the rest.
As explained in \cite[\S 8.10]{vdW70}, polynomials over $\Q$ whose Galois group
is the symmetric group $\cS_n$ can be constructed by using the fact that
the only transitive subgroup of $\cS_n$ containing a transposition and
an $(n-1)$-cycle is $\cS_n$ itself. In Proposition \ref{prop:Ken},
we force the Galois group to contain such elements
by putting together local behavior at suitable primes.
In Lemma \ref{lemma:Ken} below, we will extend
these ideas to a recursive procedure for creating a sequence
of polynomials $f_0,f_1, \ldots$ such that $\deg{f_i} = d_i$,
$\Gal{\Q(f_i)}{\Q} \simeq \cS_{d_i}$  and
$\Gal{K}{\Q} \simeq S_{d_1} \times \dots \times S_{d_n}$,
where $K$ is the compositum of the splitting fields of the $f_i$'s.
Thus the splitting field of any $f_i$ is linearly disjoint over $\Q$
from the compositum of the splitting fields of all the others.
(Actually, in Lemma \ref{lemma:Ken}, every $d_i$ will equal $7$,
but we could have used any computable sequence
$\la d_i\ra_{i\in\omega}$ instead.)

\begin{proposition}
\label{prop:Ken} 
Fix any Galois extension $E/\Q$ and any $d>1$.
Then there is a monic irreducible polynomial $f(X)$
in $\Z[X]$ of degree $d$ such that $\Gal{K}{\Q} \cong \Gal{E}{\Q} \times \cS_d$,
where $K = EF$ is the compositum of $E$ and the splitting field $F$ of $f$ over $\Q$.
In particular, $E$ and $F$ are linearly disjoint over $\Q$,
with $\Gal{F}{\Q}\cong \cS_d$.
\end{proposition}

\begin{proof}
First we recall some notation and background information.
Let $\Z_p$ be the ring of integers in the
$p$-adic field $\Q_p$ and let $\F_p$ denote the field with $p$ elements.
By Hensel's Lemma (see e.g.\ \cite[\S 18.4]{vdW70}), for a monic $h\in\Z[X]$
with mod-$(p)$ reduction $\overline{h}\in\F_p[X]$, if 
$c \in \F_p$ is a simple root of $\overline{h}$,
then there is a root $\alpha\in\Z_p$ of $h$ which,
modulo $p$, is equal to $c$ itself.  If $\overline{h}$ is
a product of distinct linear factors over $\F_p$, then $h$ splits completely
into linear factors over $\Z_p$, by applying this method to each
factor over $\F_p$.  This will be used below to satisfy the
conditions $P$ and $R$.

The finite field $\F_{q^d}$ is a Galois extension of $\F_q$ of degree $d$,
with cyclic Galois group generated by the Frobenius automorphism
$x \mapsto x^q$.  Let $\phi(X)$ be the minimal polynomial over $\F_q$
for a primitive generator of $\F_{q^d}$.  Then $\phi(X)$ has degree $d$
and splits completely in $\F_{q^d}$, with distinct roots. Now the unique
unramified extension $L$ of degree $d$ over $\Q_q$ may be constructed
as follows.  Choose $\Phi(X)\in\Z [X]$ monic of degree $d$
such that $\Phi \equiv \phi \pmod{q}$ and let $L$ be the field obtained
by adjoining a root of $\Phi$ to $\Q_q$.  Then $L$ is unramified over $\Q_q$,
since the reduction of $\Phi$ in $\F_q[X]$ is the separable polynomial $\phi$.
(In contrast, $0$ is a repeated root in $\F_q$ of the reduction of $X^2-q$,
and the splitting field $\Q_q[\sqrt{q}]$ is \emph{ramified} over $\Q_q$.)
Hensel's Lemma shows that $\Phi$ splits completely in $L$
and $\Gal{L}{\Q_q} \simeq \Gal{\F_{q^d}}{\F_q}$ also is cyclic
of order $d$.  This will be used in conditions $Q$ and $R$ below.

Now we address Proposition \ref{prop:Ken} itself.
The development here follows ideas explained more fully
in \cite[VII, \S 2]{L02}.
The Chebotarev Density Theorem \cite[\S VIII.4, Thm.\ 10]{L94}
guarantees that there are distinct primes $p,q,r \ge d$
completely split in $E/\Q$.  (This means that $E$ embeds
into each of the fields $\Q_p$, $\Q_q$, and $\Q_r$.)
Fixing these primes, we now state the conditions
we wish our polynomial $f$ to satisfy, and explain
why such an $f$ must exist.  Then we will show how
the conditions imply the theorem.

\begin{enumerate}
\item[$P$:\hspace{2 pt}] $f \equiv (X^2-\eta) u(X) \bmod{p}$,
for some $\eta\in\Z$ such that $\overline{\eta}$ is not a square
in $\F_p$ and some $u(X)\in\Z[X]$ of degree $(d-2)$ such that
$\overline{u}(X)$ splits completely into distinct linear factors over $\F_p$.
\item[$Q$:\hspace{2 pt}] $f$ is congruent modulo $q$ to the minimal
polynomial of a generator for the unique unramified extension
of degree $d$ over $\Q_q$.

\item[$R$:\hspace{2 pt}] $f \equiv X\cdot w(X) \bmod{r}$, where
$w(X)\neq X$ and $w$  is the
minimal polynomial of a generator for the unramified extension
of degree $d-1$ over $\Q_r$.
\end{enumerate}
Each condition requires that $f$ be congruent to
a particular monic polynomial of degree $d$, modulo one of the
distinct primes $p$, $q$, and $r$.
So the Chinese Remainder Theorem allows us to choose coefficients
for a monic polynomial $f \in \Z[X]$ of degree $d$
satisfying all three of these conditions.  Let $F$ be its splitting field
over $\Q$, and set $K=EF$, so both $\Gal{F}{\Q}$
and $\Gal{K}{E}$ may be seen as subgroups of $\cS_d$.
Each of the three conditions will yield a specific element of
$\Gal{K}{E}$, and the three elements together will imply
that $\Gal{K}{E}$ is all of $\cS_d$.  The process is stated in
\cite[Thm.\ VII.2.9]{L02}, and also in Example 7 of the preceding
chapter (p.\ 274).  Here we sketch it for the specific case of condition $P$.

The condition $P$ will yield a transposition
in the Galois group $\Gal{F}{\Q}$.
The polynomial $f(X)\in\Z[X]$ of degree $d$ reduces modulo $p$
to $\fbar(X)$ of the form given in condition $P$,
and the factorization there (along with the fact
that $\eta$ is not a square modulo $p$) shows that the splitting
field of $\fbar(X)$ over $\F_p$ must be a copy of $\F_{p^2}$,
with Galois group generated by the automorphism $\Phi$ of $\F_{p^2}$
which interchanges the two square roots of $\overline{\eta}$
and fixes each of the other roots of $\fbar$.
Now $X^2-\eta$ can be viewed as a polynomial
in $\Q_p[X]$, since $\Z\subset\Z_p\subset\Q_p$,
and its splitting field $L$ over $\Q_p$ is the unique
unramified extension of $\Q_p$ of degree $2$,
and in fact is the splitting field of $f$ over $\Q_p$,
since the roots of $f$ in $\F_p$ yield distinct roots
of $f$ in $\Z_p$ by Hensel's Lemma, as argued above.
The Galois group
$\Gal{L}{\Q_p} = \langle \Phi_p \rangle$
is cyclic of order $2$, with $\Phi_p$
being the lift of $\Phi$ from $\F_{p^2}$ to $L$.
This $\Phi_p$ must be a transposition, since it has order $2$ and 
must fix every other root of $f$.
For details, see \cite[Proposition VII.2.8]{L02}.
We argue next that this transposition
lifts to a transposition in $\Gal{F}{\Q}$,
by setting $e=2$ and $g(X)=X^2-\eta$ in the following lemma.

\begin{lemma}
\label{lemma:emod}
Let $E/\Q$ be Galois, $f(X)$ monic and irreducible in $\Z[X]$,
$F$ the splitting field of $f$ over $\Q$, and $K=EF$.
Assume further that $p$ is a prime completely split in $E/\Q$,
that $e>1$ is an integer, and that $g(X)\in\Z[X]$
is the minimal polynomial of a generator for the
unramified extension of degree $e$ over $\Q_p$.
If $\fbar(X)$ is the product of $\gbar(X)$ and distinct linear factors in $\F_p[X]$,
then $\Gal{K}{E}$ contains an automorphism which cyclically
permutes $e$ of the roots of $f$ and fixes each remaining root.
\end{lemma}
\begin{proof}
We can view $g(X)$ as a polynomial
in $\Q_p[X]$, since $\Z\subset\Z_p\subset\Q_p$,
and its splitting field $L$ over $\Q_p$ is the unique
unramified extension of $\Q_p$ of degree $e$.
The Galois group
$\Gal{L}{\Q_p}$ is cyclic of order $e$, generated by some
$\Phi_p$ which cyclically permutes the roots of $g$.
(Again we refer the reader to \cite[Proposition VII.2.8]{L02} for details.)
Now $F$ is the splitting field over $\Q$ of $f(X)$,
and $f(X)\in\Q_p[X]$ via the inclusion $\Z\subset\Q_p$.
Also, by assumption $E$ embeds into $\Q_p$.
We can extend this embedding to an embedding of $K$
into $L$ by noting that
$K=EF$ is generated over $E$ by the roots of $f(X)$,
which are all either roots of $g(X)$ or elements of $\Q_p$,
since by Hensel's Lemma $f(X)$ is the product of $g(X)$
with linear factors in $\Z_p[X]$.
Now the map $\Gal{L}{\Q_p}\to\Gal{K}{E}$ by restriction
(to the image of $K$ within $L$) is an injective group homomorphism,
since $L$ is generated over $\Q_p$ by the roots of $f$
in $L$.  So the restriction of $\Phi_p$ to $K$ is
the desired element in $\Gal{K}{E}$.
\qed\end{proof}

With $e=2$ and $g(X)=X^2-\eta$,
this lemma, gives us the map $\Phi_p\res K$ in $\Gal{K}{E}$
and shows it to be a transposition, as required.  
We next use Lemma \ref{lemma:emod} to satisfy conditions $Q$ and $R$.
For $Q$, we set $e=d$ and let $g(X)$ be the polynomial shown in
condition $Q$ to be congruent modulo $q$ to $f(X)$.  The lemma
(with $q$ in place of the prime $p$) shows that $\Gal{K}{E}$ contains
a cyclic permutation of order $d$, and therefore
must act transitively on the $d$ roots of $f(X)$ in $K$.
(This proves that $f(X)$ is irreducible in $E[X]$.)

Finally we take $g(X)$ to be the polynomial $w(X)$
from condition $R$, with $e=d-1$ and with $r$ as the prime
in Lemma \ref{lemma:emod}.  The lemma returns an element of
$\Gal{K}{E}$ which permutes $(d-1)$ of the roots of $f$
cyclically and fixes the last root.  But the existence of such
a permutation, along with the transposition supplied by
condition $P$ and the transitivity of the group's action,
prove that $\Gal{K}{E}\cong\cS_d$, the symmetric group
on the $d$ roots of $f(X)$ in $K$.  (See the Theorem 
on p.\ 199 in \cite[\S 8.10]{vdW70}.)
Therefore $[EF:E]=d!\geq [F:\Q]$,
making $E$ and $F$ linearly disjoint over $\Q$ and forcing
$\Gal{K}{\Q} \cong \Gal{E}{\Q} \times \cS_d$ as desired.
\qed\end{proof}

Let $F$ be any computable algebraic field.  That is, $F$ is
an algebraic field extension of its prime subfield $Q$.
The field $Q$ is isomorphic to either the $p$-element field $\F_p$
if $p=\text{char}(F)>0$, or else to the field of rational numbers.
Every one of these possibilities for $Q$ has a splitting algorithm,
and since $F$ is algebraic, this fact forces $Q$ to be computable
within $F$.  (An element $x\in F$ lies in $Q$ if and only if its minimal
polynomial over $Q$ is linear.)

Officially the domain of $F$ is $\omega$, but since the language
of fields contains the symbols $0$ and $1$ already, we will
instead write $x_0,x_1,x_2,\ldots$ for the elements of $F$.
We view $F$ as the union of an infinite chain of finitely generated subfields:
$$ Q=F_0 \subseteq F_1\subseteq F_2\subseteq\cdots \subseteq F,$$
where $F_s = Q(x_0,\ldots,x_{s-1})$ for every $s$.
The Effective Theorem of the Primitive Element (see \cite{FJ86},
or \cite[Theorem 3.11]{MS11}) allows
us to compute for each $s$ a single element $z_s\in F_s$
that generates all of $F_s$; indeed we may assume $z_s=x_t$
for the least $t$ such that $F_s=Q(x_t)$.
We may compute the minimal polynomial $q_s(X)\in Q[X]$
of each $z_s$ over $Q$, and also compute polynomials
$p_s\in Q[X_0,\ldots,X_{s-1}]$ such that $p_s(z_0,\ldots,z_{s-1},X)$
is the minimal polynomial of $z_s$ over $F_{s-1}$.

Now, following \cite{M09b} and \cite{MS11}, we define the
\emph{automorphism tree for $F$} to be the following subtree of $\omega^{<\omega}$:
$$ I_F = \set{\sigma\in \omega^{<\omega}}{
(\forall s<|\sigma|)~p_{s}(x_{\sigma(0)},\ldots,x_{\sigma(s)})=0}.$$
More generally, for a computable field $E$ isomorphic to $F$,
with domain $\{ y_0,y_1,\ldots\}$,
we define the \emph{isomorphism tree} for $F$ and $E$ to be
$$ I_{FE} = \set{\sigma\in \omega^{<\omega}}{
(\forall s<|\sigma|)~\ptilde_{s}(y_{\sigma(0)},\ldots,y_{\sigma(s)})=0},$$
where $\ptilde_s\in\Qtilde [X_0,\ldots,X_s]$
is the image of $p_s$ when its coefficients are mapped
into the prime subfield $\Qtilde$ of $E$ by the (unique)
isomorphism between these prime subfields.
It follows that isomorphisms from $F$ onto $E$
correspond precisely to paths through $I_{FE}$,
and that this correspondence preserves Turing degrees.
(For an isomorphism $h$, the path contains those nodes $\sigma\in I_{FE}$
with $y_{\sigma(i)}=h(z_i)$ for all $i<|\sigma|$; conversely, given a path
$f$ through $I_{FE}$, define $h(z_s)=y_{f(s)}$ for all $s$.)
All of this material is described in detail in \cite{M09b} and 
\cite{MS11}.

The isomorphism tree $I_{FE}$ can be defined the same way
(but perhaps should be renamed) for any algebraic
fields $E$ and $F$, whether or not they are isomorphic.
In this more general case, paths through $I_{FE}$
correspond to field embeddings of $F$ into $E$.
If no such embeddings exist, then $I_{FE}$ must be a finite tree.
The following useful corollary generalizes this statement.

\begin{lemma}
\label{lemma:finiteblock}
Let $F$ and $E$ be algebraic fields of the same characteristic,
and suppose that $x\in F$ and $y\in E$ have the same minimal
polynomial over the prime subfields.
If there is no field embedding $F\to E$ with $x\mapsto y$,
then there exists some finitely generated subfield $F_0\subseteq F$
containing $x$ such that no field embedding of $F_0$
into $E$ maps $x$ to $y$.
\end{lemma}
\begin{proof}
Build the tree $I_{FE}$ as above, using an
enumeration of $F$ that lists $x$ as its first element.
Then $\la y\ra$ constitutes a node on this tree, through
which there is no path.  K\"onig's Lemma shows that
$I_{FE}$ contains only finitely many nodes above $\la y\ra$,
all of height $<n$, say.  Therefore, the field $F_0$
generated by the first $n$ elements of $F$ (in our enumeration)
has no embedding into $E$ with $x\mapsto y$.
\qed\end{proof}
\begin{corollary}
\label{cor:nofiniteblock}
Two algebraic fields $E$ and $F$ are isomorphic if and only if every
finitely generated subfield of each one embeds into the other.
\end{corollary}
\begin{proof}
Apply Lemma \ref{lemma:finiteblock}, with $x$ and $y$ being the zero
elements of their fields, to see that each field
embeds into the other.  Then apply Lemma \ref{lemma:endo}
to the composition of the two embeddings.
\qed\end{proof}

\begin{corollary}
\label{cor:iso}
Two algebraic fields $E$ and $F$ are isomorphic if and only if
they have the same characteristic and every polynomial
over the prime subfield with a root in either $E$ or $F$
also has a root in the other.  (Here we extend the unique
isomorphism between the prime subfields to an isomorphism
of the polynomial rings over those subfields.)
\end{corollary}
\begin{proof}
For the nontrivial direction, we apply the Theorem of the
Primitive Element, which states that each finitely
generated subfield (of $E$, say) is generated by a single
element $x$.  Let $p(X)$ be the minimal polynomial of $x$
over the prime subfield $Q_E$ of $E$.  The isomorphism
$f:Q_E\to Q_F$ maps $p(X)$ to an irreducible polynomial
$p^f(X)\in Q_F[X]$, which by hypothesis must have a root $y\in F$.
But then $Q_E(x)\cong Q_E[X]/(p(X))\cong Q_F[X]/(p^f(X))\cong Q_F(y)$.
So every finitely generated subfield of $E$ embeds into $F$,
and conversely.  Corollary \ref{cor:nofiniteblock} completes the proof.
\qed\end{proof}

Finally, we must also give some technical definitions for embeddings
among fragments of fields, by which we mean subsets of a field
that are not necessarily subfields themselves.  This notion will arise 
when we examine all computable copies of a given field:  we will need
to consider the structure (if any) computed by each partial
computable function $\phi_e$, without knowing for certain whether
it is a field or not.

\begin{definition}
\label{defn:frag}
Let $Q$ be the prime field of a given characteristic,
and $\overline{Q}$ a computable presentation of its algebraic closure,
with $Q$ regarded as a subfield.
Suppose that $C$ is a subset of $\omega$ on which two binary operations,
which we denote by $+$ and $\cdot$, are partially defined, i.e., for $x,y\in C$, the
values of $x+y$ and $x\cdot y$ may be other elements of $C$, or may be undefined.
(These operations are not the usual addition and multiplication on $\omega$;
they are simply any binary operations, and we wish to determine
whether the set $\omega$ forms a field under these operations.)
Assume that $C$ contains distinct elements $x$ and $y$
satisfying $x+x=x$ and $y\cdot y=y$; we refer to $x$ as $0$
and to $y$ as $1$, since they will represent identity elements
in this fragment of a field.  If these partial operations
do not contradict associativity, commutativity, distributivity,
cancellation, the identity properties of $0$ and $1$, or equality
of the characteristics of $C$ and $Q$, then $C$ is a \emph{field fragment}
of that characteristic.  To contradict associativity, for example,
would require the existence of three elements $c_0,c_1,c_2\in C$
such that the four products $(c_0\cdot c_1)$, $(c_1\cdot c_2)$,
$(c_0\cdot c_1)\cdot c_2$, and $c_0\cdot (c_1\cdot c_2)$ are all defined,
but the last two of these products are distinct.  Similarly, a contradiction
to cancellation for $\cdot$ would consist of elements $c_0\neq 0$ and
$c_1\neq c_2$ in $C$ with $(c_0\cdot c_1)$ and $(c_0\cdot c_2)$ 
both defined and equal to each other.  All the other properties mentioned
are also universal statements (given that $0$ and $1$ have already been specified),
and so a contradiction to any of them is simply an instantiation of the negation.
(Notice that we do not check for the existence of inverses.)

For a field fragment $C$,
we define the prime field fragment $Q_0\subseteq C$ as follows.
If $Q$ has positive characteristic, then $Q_0$ contains $0$
and all elements $1,(1+1),(1+1+1),\ldots$ that are defined in $C$.
If $Q$ has characteristic $0$, fix some enumeration $\{ q_0,q_1,\ldots\}$
of the rationals $\Q$, For each $q_i=\frac{m}{n}$ in turn,
we check whether $C$ possesses an element of the form
$\frac{1+1+\cdots+1}{1+\cdots +1}$, with $m$ $1$'s on top
and $n$ $1$'s below, and also possesses an additive inverse of
that element.  If so, we enumerate both elements into $Q_0$.
If not, then we stop enumerating entirely, so $Q_0$ consists of only
those elements enumerated so far.  We make the obvious identification
between elements of $Q_0$ and elements of $\Q$, noting 
that this identification must be unique, since $C$ is a field fragment.

Suppose that $D=\{ x_0,\ldots, x_m\}\subseteq C$.
An \emph{embedding of $D$ into $\overline{Q}$} is a function $f:D\to\overline{Q}$
such that for every $i\leq m$, there is some $p_i\in Q_0[X_0,\ldots,X_i]$
with $p_i(x_0,\ldots,x_i)=0$ in $C$
(specifically, all sums and products in this calculation lie in $C$),
such that $\overline{p_i}(f(x_0),\ldots,f(x_{i-1}),X_n)$
is irreducible in $Q(x_0,\ldots x_{i-1})[X_n]$ and
$\overline{p_i}(f(x_0),\ldots,f(x_i))=0$ in $\overline{Q}$,
where $\overline{p_i}\in\overline{Q}[X_0,\ldots,X_i]$ is the image of $p_i$
via the unique embedding of $Q_0$ into $\overline{Q}$.

Finally, if $C'$ is another field fragment and $D'\subseteq C'$ is finite,
then an \emph{embedding of $D'$ into $D$} consists of embeddings
$f$ of $D$ into $\overline{Q}$ and $f'$ of $D'$ into $\overline{Q}$ such that
$\rg{f'}\subseteq\rg{f}$.  (We often think of $f^{-1}\circ f'$ as the actual embedding.)
\end{definition}

The resulting lemma is easily proven.

\begin{lemma}
\label{lemma:frag}
For finite subsets $D$ and $D'$ of finite field fragments $C$ and $C'$,
it is decidable whether there exists an embedding of $D'$ into $D$.
The procedure is uniform in the (partial) operations on $C$ and $C'$, but 
it is necessary to know the exact size of each finite set $C$, $C'$, $D$, and $D'$,
and to be able to decide the domains of the partial operations.
\end{lemma}

\section{Relative Computable Categoricity}
\label{sec:rcc}

The definition of a computably categorical structure is often
strengthened to consider
noncomputable copies, yielding the notion of relative computable
categoricity introduced in Definition \ref{defn:cc}.
Relative computable categoricity is actually a more natural property than
computable categoricity, in the sense that, by results in \cite{AKMS89}
and \cite{Chis}, it has a fairly simple syntactic characterization:
a structure is relatively computably categorical
if and only if it has a $\Sigma^0_1$ Scott family over finitely many
parameters.

\begin{definition}
A \emph{$\Sigma^0_1$ Scott family} for a structure $\A$ is a
computable sequence $\theta_0(\bar a, x_0, \ldots, x_{n_0}),
\theta_1(\bar a, x_0, \ldots, x_{n_1}), \ldots$ of $\exists$-formulas,
where $\bar a$ is a finite tuple from $\A$, such that every tuple of
elements from $\A$ satisfies at least one $\theta_i$, and for each
$i$, any two tuples satisfying $\theta_i$ can be interchanged by an
automorphism of $\A$.
\end{definition}

In contrast, no such nice characterization of computable categoricity
is known. Indeed, while the property of having a $\Sigma^0_1$ Scott family
is arithmetically a $\Sigma^0_3$ property, computable categoricity
is known by work in \cite{W03} to be at least $\Pi^0_4$-hard.
(In fact, we will show in Section \ref{sec:cc} that it is this hard
just within the class of algebraic fields.)  In this section and the next
two, we add further evidence in favor of the neatness of relative
computable categoricity, by considering computable algebraic fields.

Following the notation of \cite{MS11}, we denote the \emph{full orbit relation}
on $F$ by
$$A_F = \set{\la a_1,\ldots,a_n,b_1,\ldots,b_n\ra\in F^{<\omega}}{
(\exists\sigma\in\Aut{F})(\forall i\leq n)~\sigma(a_i)=b_i}.$$
The simple \emph{orbit relation} then is $B_F = A_F\cap F^2$,
the binary relation on $F$ of being in the same orbit.
By Lemma \ref{lemma:finiteblock}, if a pair $\la x,x'\ra$
does not lie in $B_F$, then there is some finitely generated subfield
$F_t\subseteq F$ containing $x$ such that no embedding
of $F_t$ into $F$ can map $x$ to $x'$.  (Here
$Q=F_0\subseteq F_1\subseteq\cdots\subseteq F$ can be
any chain of subfields with union $F$.  Normally we will
set $F_s$ to be the subfield generated by the first $(s-1)$
elements of $F$.)  Of course, if $x$ and $x'$ are
not conjugate over $Q$, then we will realize immediately
that $\la x,x'\ra\notin B_F$; the interesting case is that in
which $x$ and $x'$ are conjugate, yet lie in distinct orbits.
Since there are only finitely many conjugates of $x$ in $F$,
the following definition makes sense.

\begin{definition}
\label{defn:orbitfct}
Let $F$ be a computable algebraic field, with prime subfield $Q$
and such that $F$ is the union of a chain
$Q=F_0\subseteq F_1\subseteq\cdots\subseteq F$ of finitely generated subfields.
The \emph{orbit function} $h:F\to\omega$ of $F$ defines $h(x)$ to equal
the least $t$ such that $x\in F_t$ but, for every $x'\in F$ that is 
conjugate to $x$ over
$Q$ with $\la x,x'\ra\notin B_F$, there is no embedding of $F_t$ into $F$
that maps $x$ to $x'$.
\end{definition}

So the function $h$, on input $x$, considers all false conjugates $x'$
of $x$ in $F$:  that is, all conjugates of $x$ that do not
lie in the orbit of $x$ under automorphisms.  (This notion is important
for categoricity, since the false conjugates in $F$ correspond to
elements of computable copies of $F$ to which we might mistakenly
map $x$.)  Therefore, $h$ has the property that, for all $x\in F$,
$$ [x\mapsto x' \text{~extends to~}F_{h(x)}] \text{~~if and only if~~}
\la x,x'\ra\in B_F.$$
Moreover, $h(x)$ is the least number with this property.
Of course, $x$ might not have any false conjugates in $F$;
in this case $h(x)$ is the least $s$ with $x\in F_s$.

We note that the function $h$ does depend on
the choice of chain $Q=F_0\subseteq F_1\subseteq\cdots\subseteq F$.
Normally we take this to be a computable chain
(that is, we can compute a strong index for a finite generating set
for each $F_s$, uniformly in $s$), and it is readily proven
that with the orbit function for one computable chain as an oracle,
we can compute an upper bound for the orbit function for any other
computable chain.  However, for many fields the orbit function
will not be computable.  Indeed, should it be computable,
then $F$ must be computably categorical.
More generally, we have the following lemma.
\begin{lemma}
\label{lemma:rcc}
For a computable algebraic field $F$, if the orbit function $h$
(for some computable chain $Q=F_0\subseteq F_1\subseteq\cdots$
of finitely generated subfields with union $F$)
is computably bounded, then $F$ is relatively computably categorical.
\end{lemma}
\begin{proof}
Fix some computable function $g:F\to\omega$
such that $g(x)\converges\geq h(x)$ for all $x$.
It is straightforward to construct a $\Sigma^0_1$ Scott family
for $F$, fixing a primitive generator $z_s$ for each $F_s$ and
its minimal polynomial $q_s(Z)\in Q[Z]$.  We also fix,
for each $s<t$, a polynomial $p_{s,t}\in Q[Z,Y]$ such that
$p_{s,t}(z_s,Y)$ is the minimal polynomial of $z_t$
over the subfield $F_s=Q(z_s)$.  (This definition can be done 
effectively.)

For each $z_s$, we have a formula $\gamma_s(Z)$ saying:
$$ q_s(Z)=0~\&~\exists y~p_{s,g(s)}(Z,y)=0.$$
By Definition \ref{defn:orbitfct}, and since $g(z_s)\geq h(z_s)$,
every $z$ satisfying $\gamma_s(z)$ lies in the orbit of $z_s$.
For a general tuple of elements $x_0,\ldots,x_n\in F$,
we find the least $s$ with all $x_i\in F_s=Q(z_s)$, find polynomials
$r_i(Z)\in Q[Z]$ with $r_i(z_s)=x_i$, and define the formula
$\delta_{\xvec}(X_0,\ldots,X_n)$:
$$ \exists z [\gamma_s(z)~\&~(\forall i\leq n)X_i=r_i(z)].$$
The Scott family is then just the set $\Ss=\set{\delta_{\xvec}(\Xvec)}{\xvec\in F^{<\omega}}$,
In fact, it is a particularly nice Scott family,
since it involves no parameters from $F$. (Of course,
polynomial equations with parameters from $Q$ can be expressed
entirely in terms of the constant symbols $0$ and $1$ and the
operations of addition and multiplication.)

Since relative computable categoricity is equivalent to having
a $\Sigma^0_1$ Scott family, this definition completes the proof of
Lemma \ref{lemma:rcc}.  However, to illuminate the use of
the orbit function further, we will also give a direct
proof of relative computable categoricity of $F$, by constructing
a computable isomorphism from $F$ onto an arbitrary computable field
$\Ftilde$ isomorphic to $F$.  Our construction relativizes easily
to the degree of a noncomputable field $\Ftilde$.
Of course, the prime subfield $Q=F_0$ is computable within
the algebraic field $F$ and has a unique embedding $f_0$
into $\Ftilde$; this $f_0$ is computable and has image $\Qtilde$,
the prime subfield of $\Ftilde$.  Moreover, $f_0$ is known
to extend to an isomorphism from $F$ onto $\Ftilde$,
since the two fields are assumed to be isomorphic
and every isomorphism between them must restrict to $f_0$.

Now assume inductively that we have constructed an embedding
$f_s:F_s\to\Ftilde$ that extends to some (not necessarily computable)
embedding $\rho$ of $F$ into $\Ftilde$.
Let $\qtilde_{s+1}$ and $\ptilde_{s+1,g(s+1)}$ be the images
of the polynomials $q_{s+1}$ and $p_{s+1,g(s+1)}$ under the
map $f_0:Q\to\Qtilde$ on their coefficients.  To extend $f_s$ to $F_{s+1}$,
compute $g(s+1)$ and search for any two elements $\ytilde,\ztilde\in\Ftilde$
such that $\qtilde_{s+1}(\ztilde)=0$ and
$\ptilde_{s+1,g(s+1)}(\ztilde,\ytilde)=0$,
and such that the map $z_{s+1}\mapsto\ztilde$ would send $z_s$ to $f_s(z_s)$.
(Of course, $z_s\in F_{s+1}=Q(z_{s+1})$, so the choice of an image
for $z_{s+1}$ uniquely determines the image of $z_s$.  Also,
notice that these conditions are satisfied when
$\ytilde=\rho(z_{h(s+1)})$ and $\ztilde=\rho(z_{s+1})$,
so we must eventually find some $\ytilde$ and $\ztilde$ as desired,
although they will not necessarily be $\rho(z_{h(s+1)})$ and $\rho(z_{s+1})$.)
Define $f_{s+1}(z_s)$ to be this $\ztilde$, thus defining $f_{s+1}$
on all of $F_{s+1}$.  We claim that this $f_{s+1}$ extends $f_s$
and is a field embedding of $F_{s+1}$ into $\Ftilde$, and also that
$f_{s+1}$ itself extends to some embedding of all of $F$ into $\Ftilde$.

First, since $\qtilde_{s+1}(\ztilde)=0$ and $q_{s+1}(z_{s+1})=0$,
and since $q_{s+1}$ and $\qtilde_{s+1}$ are irreducible
in $Q[Z]$ and $\Qtilde[Z]$, respectively, we know that
$$ F_{s+1} = Q(z_{s+1}) \cong Q[Z]/(q_{s+1}(Z)) \cong
\Qtilde[Z]/(\qtilde_{s+1}(Z)) \cong \Qtilde(\ztilde),$$
via the map $z_{s+1}\mapsto\ztilde$.  Therefore $f_{s+1}$ really
is a field embedding of $F_{s+1}$ into $\Ftilde$.
Moreover, we checked that $z_{s+1}\mapsto\ztilde$ sends $z_s$ to $f_s(z_s)$,
and so $f_s\subseteq f_{s+1}$.  It remains to see that $f_{s+1}$
extends to a field embedding $\alpha:F\to\Ftilde$.
To prove this fact, notice that with $\ptilde_{s+1,g(s+1)}(\ztilde,\ytilde)=0$,
we must have $p_{s+1,g(s+1)}(z,y)=0$ as well, where
$z=\rho^{-1}(\ztilde)$ and $y=\rho^{-1}(\ytilde)$.
Since $g(s+1)\geq h(s+1)$, this fact means that
$p_{s+1,h(s+1)}(z,Y)$ must have a root in $F$ as well.
But by Definition \ref{defn:orbitfct}, we then have
$\la z_{s+1},z\ra\in B_F$, so there is some automorphism
$\beta$ of $F$ with $\beta(z_{s+1})=z$.  But then
$$ \rho(\beta(z_{s+1})) = \rho(z) = \ztilde,$$
and so $(\rho\circ\beta)$ is an embedding of $F$ into $\Ftilde$
that extends $f_{s+1}$, as required.

Therefore, the union $f=\bigcup_s f_s$ is a well-defined
computable function with domain $\bigcup_s F_s = F$.
Since every $f_s$ is a field embedding, so is $f$.
But with $\Ftilde\cong F$, a field embedding of $F$
into $\Ftilde$ must in fact be an isomorphism,
by Corollary \ref{cor:endo}.  This conclusion
proves Lemma \ref{lemma:rcc}.
\qed\end{proof}

We can also prove the converse of Lemma \ref{lemma:rcc},
yielding the full equivalence.
\begin{theorem}
\label{thm:rcc}
For a computable algebraic field $F$, the following are equivalent.
\begin{enumerate}
\item
$F$ is relatively computably categorical.
\item
$F$ has a  $\Sigma^0_1$ Scott family.
\item
$F$ has a $\Sigma^0_1$ Scott family
using no parameters from $F$.
\item
The orbit function $h$ for $F$ (with respect to some
computable chain $Q=F_0\subseteq F_1\subseteq\cdots$) is computably bounded.
\item
The orbit function $h$ for $F$ (with respect to every
computable chain) is computably bounded.
\end{enumerate}
\end{theorem}
\begin{proof}
$(5)\implies (4)$ is immediate, and
the equivalence of (1) and (2) was established
(for all computable structures, not just fields)
by Ash, Knight, Manasse, and Slaman in \cite{AKMS89},
and independently by Chisholm in \cite{Chis}.
Lemma \ref{lemma:rcc} shows $(4)\implies (1)$,
and its proof also explained how (4) yields (3), which in turn
clearly implies (2).  So we now prove $(2)\implies (5)$.
Fix any computable chain with union $F$, and fix
a $\Sigma^0_1$ Scott family $\Ss$ for $F$.  By the Effective Theorem
of the Primitive Element, we may assume that $\Ss$ uses
(at most) a single parameter $a$ from $F$. 
Since $a$ has finitely many conjugates over $Q$ in
the algebraic field $F$, we may assume that we know
all elements $a=b_0,b_1,\ldots,b_m$ satisfying
$\la a,b\ra\in B_F$.  Also, for each $b\in F$ conjugate
to $a$ but with $\la a,b\ra\notin B_F$, Corollary
\ref{lemma:finiteblock} yields an $s$ such that no embedding
of $F_s$ into $F$ maps $a$ to $b$.  So we may fix some $s_0$
so large that all conjugates of $a$ in $F$ lie in $F_{s_0}$,
but that $b_0,\ldots,b_m$ are the only possible images of $a$
under embeddings of $F_{s_0}$ into $F$.

We now compute the value $g(z_s)$ for a primitive generator
$z_s$ of $F_s$, for arbitrary $s > s_0$.
For each $i\leq m$, we list out the formulas of $\Ss$ and search through
$F$ for witnesses for these existential formulas, until we find a formula
$$\gamma_i(Z) =  \exists\xvec \delta_i(Z,\xvec,a) \in \Ss,$$
an $n_i$, and a tuple $\vvec_i$ from $F$ such that
$\delta_i(z_s,\vvec_i,b_i)$ holds in $F$.
It is important to notice that in our search, we have replaced
the parameter $a$ by $b_i$.  (Of course $b_0=a$, so
one search involving $a$ does still take place.)
Now for each $i$, there is (at least one) automorphism
$\psi_i$ of $F$ with $\psi_i(b_i)=a$.
By the definition of Scott family, there is a
formula $\gamma(Z)=\exists\xvec \delta(Z,\xvec,a)$ in $\Ss$
for which $\gamma(\psi_i(z_s))$ holds; that is,
$\delta(\psi_i(z_s),\vvec,a)$
holds for some $\vvec$.  But then $\delta(z_s,\psi_i^{-1}(\vvec),b_i)$ also
holds, since $\psi_i$ is an automorphism.  Therefore eventually
our search halts and produces a formula $\gamma_i(Z)$,
a tuple $\vvec_i$, and an $n_i$.
Define $g(z_s)$ to be the least number $\geq s$ such that 
for all $i\leq m$, the entire
finite tuple $\vvec_i$ is contained in the subfield $F_{g(z_s)}$.
This definition completes our computation of $g$, on every $s > s_0$;
now we must show that $g(z_s)\geq h(z_s)$.

So consider any conjugate $z$ of $z_s$,
such that $\la z_s,z\ra\notin B_F$.
Now for each $i$, $\gamma_i(\psi_i(z_s))$ holds in $F$
(with $\psi_i$ as chosen above).  But all $\psi_i(z_s)$ lie in the orbit of $z_s$,
and all $\psi_i(z)$ lie in the orbit of $z$.
Since these two orbits are distinct, the definition of Scott family
shows that $\gamma_i(\psi_j(z))$ must be false in $F$ for each $i,j\leq m$.
So there is no tuple $\xvec$ from $F$
for which any $\delta_i(\psi_j(z),\xvec,a)$ holds,
and hence (by applying the automorphism $\psi_j^{-1}$)
no tuple $\xvec$ for which any $\delta_i(z,\xvec,b_j)$ holds either.

Now if $p_{s,g(s)}(z,Y)$ had a root $y$ in $F$ (where $p_{s,g(s)}$ is
as in the proof of Lemma \ref{lemma:rcc}),
then there would be an isomorphism from $F_{g(z_s)}$
onto $F(y)$:  start with a field embedding $\psi$ of $F_s$ into $F$
that sends $z_s$ to its $F_0$-conjugate $z$, and then extend
by setting $\psi(z_{g(z_s)})=y$ (which is still a field embedding of $F_{g(z_s)}$
into $F$, since $z_{g(z_s)}$ has minimal polynomial $p_{s,g(z_s)}(z_s,Y)$
over $F_s$ and $y$ has minimal polynomial $p_{s,g(z_s)}(\psi(z_s),Y)$
over $\psi(F_s)$).  However, then the quantifier-free formula
$$ \delta_0 (z,\psi(\vvec_0),\psi(a))$$
must hold in $F$, since $\delta_0(z_s,\vvec_0,a)$ held
in $F_{g(z_s)}$.  By our choice of $s_0$, and
since $s\geq s_0$, we have $\psi(a)=b_j$ for some $j$,
contradicting our conclusion above
that $\delta_i(z,\xvec,b_i)$ fails for all tuples $\xvec$ in $F$.

So, for all conjugates $z$ of $z_s$ with
$\la z_s,z\ra\notin B_F$, the polynomial $p_{s,g(z_s)}(z,Y)$ has no root in $F$.
By the minimality of $h(z_s)$ in Definition \ref{defn:orbitfct},
this fact forces $g(z_s)\geq h(z_s)$.  The argument above
only defined $g$ on inputs $z_s$ with $s > s_0$, but it requires just finitely
much more information to set $g(z_s)=h(z_s)$ for all $s\leq s_0$.
Then one simply defines $g(x)$ on arbitrary $x\in F$
by $g(x)=g(z_s)$, where $s$ is minimal such that $x\in F_s$.
This definition works because, if $F_{g(x)}$ has an embedding $\alpha$ into $F$
with $x'=\alpha(x)$, then by our construction,
$\la z_s,\alpha(z_s)\ra\in B_F$.  Hence the map $\alpha\res F_s$ extends
to an automorphism of $F$, and since $x\in F_s$,
that automorphism maps $x$ to $\alpha(x)=x'$.
So the computable function $g$ does indeed bound $h$.
\qed\end{proof}

As mentioned in the introduction, it was shown in \cite{MS11} that if
the algebraic field $F$ has a splitting algorithm, then $F$ is
(relatively) computably categorical if and only if $B_F$ is computable. We now
show that the situation is different for a general algebraic field $F$.
Using Scott families, we easily show that relative computable
categoricity implies that $B_F$ is $\Sigma_1$.
The challenging part is to show that the converse is false,
and here we actually strengthen that result by making $B_F$
(and also the full orbit relation $A_F$) computable.
\begin{theorem}
\label{thm:rccorbit}
Let $\M$ be any computable, relatively computably categorical structure.
Then the full orbit relation of $\M$ (i.e., the set of pairs of tuples
$\la \avec,\bvec\ra$ such that some automorphism of $\M$ maps $\avec$
to $\bvec$) is computably enumerable.
However, there exists a computable algebraic field $F$
that is not computably categorical, yet has $A_F$ computable.
\end{theorem}
\begin{proof}
Suppose that $\M$ is a computable, relatively computably categorical
structure, and denote its full orbit relation by $A_{\M}$.
Let $\Ss$ be a $\Sigma^0_1$ Scott family for $\M$.  For any pair
$\la \avec,\bvec\ra$ of $n$-tuples from $\M$,
search through all formulas in $\Ss$ with exactly
$n$ free variables, and all tuples of possible witness
elements for each formula.  Enumerate
$\la\avec,\bvec\ra$ into $A_{\M}$ if ever we find a single
formula in $\Ss$ satisfied by both these tuples.

Of course, every $\avec$ satisfies some formula in $\Ss$,
and if $\la\avec,\bvec\ra\in A_{\M}$, then $\bvec$ satisfies that
same formula, so the pair is enumerated.  Conversely,
by the definition of Scott family, any two tuples satisfying
the same formula must be images of each other
under automorphisms of $\M$.

We describe a simple version of the basic module used
to construct the field $F$ and its computable copy
$\Ftilde$ which together prove the second statement.
To ensure that a single $\phi_e$ is not an isomorphism from
$F$ onto $\Ftilde$, we use the cube roots of $2$.
Each field starts with one $\sqrt[3]2$, called $\theta_0$ and
$\tilde{\theta}_0$, respectively.
In order to be an isomorphism, $\phi_e(\theta_0)$
must converge to $\tilde{\theta}_0$.  If this convergence ever happens,
then we adjoin two more cube roots $\theta_1$ and $\theta_2$
of $2$ to $F$, and likewise in $\Ftilde$.
In $F$ we also \emph{tag} the original $\theta_0$,
by finding a polynomial $q\in\Q[X,Y]$ such that we can
adjoin a root of $q(\theta_0,Y)$ to $F$ without adjoining
any roots of either $q(\theta_1,Y)$ or $q(\theta_2,Y)$;
thus these two conjugates of $\theta_0$
are not tagged.  In $\Ftilde$ we adjoin a root of $q(\tilde{\theta}_1,Y)$.
Thus the two fields remain isomorphic, but only via isomorphisms
mapping $\theta_0$ to $\tilde{\theta}_1$.

Moreover, no matter what the outcome of this basic module, 
every computable field $E\cong F$ has computable orbit relation $B_E$.
Certainly every pair $\la x,x\ra\in B_E$,
and if our program for computing $B_E$ is ever given a
pair $\la\theta,\theta'\ra$ of distinct cube roots of $2$, then it searches for
a third such cube root $\theta''$ and also for a tag,
i.e.\ a root of $q(\theta,Y)$ or of $q(\theta',Y)$ or of $q(\theta'',Y)$.
These must exist, since the existence of the second cube root
means that the basic module must have performed the action
described above.  Once our program finds the tag, it knows
which two of these roots lie in the same orbit, and therefore
can answer correctly whether $\la\theta,\theta'\ra\in B_E$.

In fact, such tagging requires an elaborate algebraic proof
of the existence of the appropriate polynomials.
(For instance, letting $q(X,Y)$ be $(Y^2-X)$
would not work, because any square root of $\theta_0$
would generate square roots of $\theta_1$ and $\theta_2$ as well.)
Moreover, if any basic module adjoins to $F$ all three cube roots of some element
(such as $2$), then every $x$ with a cube root $F$ must
have three cube roots in $F$.  Thus this basic module
does not extend readily to infinitely many requirements,
so we resort to a similar strategy using roots of polynomials
with symmetric Galois groups $\cS_7$ of order $7$.

Viewing such a Galois group as the symmetric
permutation group on the seven roots of some polynomial,
we see that its symmetry allows us to adjoin the sum of any subset $S$
(with $1<|S|<6$) of the seven roots we like while keeping
the individual roots out of the field.  For instance,
if $x$ and $y$ are two of these roots, then the subfield
generated by $(x+y)$ is the fixed field of the subgroup
of $\cS_7$ containing those permutations fixing $\{ x,y\}$
setwise, and as this subgroup fixes no individual root,
the subfield contains none of the roots.
Moreover, $\cS_7$ is a sufficiently
large group to allow us both to adjoin specific elements
and to tag them.  By analogy to the basic module above,
think of $\theta_0=x+y$ as the sum of two roots $x$ and $y$,
and $\theta_1=u+v$ as the sum of two others $u$ and $v$.
We can tag either one later if necessary by adjoining $x$ or $u$.
If there were only four roots in total, then adjoining $x+y$
would have forced $u+v$ to enter the field as well.
With only five, adjoining both $(x+y)$ and $(u+v)$
would force the fifth root to enter the field, since every permutation
fixing both $\{ x,y\}$ and $\{ u,v\}$ would have to fix
the fifth root.  With only six,
it would force the sum of the other two roots to enter the field,
where it would be conjugate to $\theta_0$ and $\theta_1$,
in the manner of $\theta_2$ above; this would not ruin
the argument, but it would complicate it, so we use
seven roots instead.

\begin{lemma}
\label{lemma:Ken}
There exists a computable infinite sequence $p_0(X),p_1(X),\ldots$
of polynomials of degree $7$ in $\Q[X]$ such that for every $e\in\omega$,
if $K_e$ is the subfield of $\Qbar$ generated by the splitting fields
of all $p_i(X)$ with $i\neq e$, then the splitting field $P_e$ of $p_e(X)$
over $K_e$ is the symmetric group on the seven roots of $p_e(X)$.
\end{lemma}
\begin{proof}
Apply Proposition \ref{prop:Ken} repeatedly,
to get $\Q=E_0\subset E_1\subset E_2\subset\cdots$ and
polynomials $p_0,p_1,\ldots$, all of degree $7$ in $\Z[X]$,
such that each $E_{e+1}$ is the compositum of $E_e$
with the splitting field $P_e$ of $p_e$ over $\Q$
and such that $\Gal{P_e}{\Q}\cong\cS_7$ for all $e$.
Notice that if $P_e$ were not linearly disjoint from $K_e$
(as defined in Definition \ref{defn:lindisjt}), then for some $j$, 
the field $P_e$ 
would fail
to be linearly disjoint from the subfield $K_e\cap E_{j+1}$,
which (for the least such $j$) would contradict
the linear disjointness of $P_j$ from $E_j$.
\qed\end{proof}

We fix a sequence of polynomials $p_e(X)$ as described
in Lemma \ref{lemma:Ken}.  In the following construction,
using a fixed computable copy $\Qbar$ of the algebraic closure of $\Q$,
we let $x_e$, $y_e$, $u_e$, and $v_e$ be the four $<$-least roots of
$p_e$ in $\Qbar$.  Our fields $F$ and $\Ftilde$ will both be computably
enumerable subfields of $\Qbar$, hence computably isomorphic
(by taking pullbacks) to computable algebraic fields of characteristic $0$.
The requirements, which never injure one another, are
$$\R_e:~~\phi_e\text{~is not an isomorphism from $F$ onto~}\Ftilde,$$
for every $e\in\omega$.  

With no injury involved, we may define our fields quite simply.
The c.e.\ subfield $F$ contains precisely the following elements of $\Qbar$:
\begin{itemize}
\item
$(x_e+y_e)$, for all $e\in\omega$;
\item
$(u_e+v_e)$, for all $e\in\omega$ for which $\phi_e(x_e+y_e)\converges=x_e+y_e$;
\item
$x_e$ (hence also $y_e$),
for all $e\in\omega$ for which $\phi_e(x_e+y_e)\converges=x_e+y_e$;
\item
and all elements of $\Qbar$ generated by these.
\end{itemize}
$\Ftilde$ contains precisely the following elements of $\Qbar$:
\begin{itemize}
\item
$(x_e+y_e)$, for all $e\in\omega$;
\item
$u_e$ and $v_e$ (hence also $(u_e+v_e)$),
for all $e\in\omega$ with $\phi_e(x_e+y_e)\converges=x_e+y_e$;
\item
and all elements of $\Qbar$ generated by these.
\end{itemize}

We readily define an isomorphism $\rho$ from $F$ onto $\Ftilde$.
If $\phi_e(x_e+y_e)\converges = x_e+y_e$, then
$\rho(x_e)=u_e$, $\rho(y_e)=v_e$, and $\rho(u_e+v_e)=x_e+y_e$.
Otherwise $\rho(x_e+y_e)=x_e+y_e$.
(In fact, there are $2^\omega$-many isomorphisms, since
the first case occurs for infinitely many $e$, and when it does,
$\rho(x_e)=v_e$ and $\rho(y_e)=u_e$ is also possible; moreover,
by Lemma \ref{lemma:Ken}, the choices of $\rho(x_e)$ may be made
independently for all the different $e$ for which the first case occurs.)

Let $F_e\subseteq F$ be generated by $\{ (x_e+y_e),(u_e+v_e),x_e\}\cap F$,
and $\Ftilde_e\subseteq\Ftilde$ by $\{(x_e+y_e),u_e,v_e\}\cap\Ftilde$.
Let $L_e =  K_e \cap F$ be the subfield of $F$ generated by all
those generators of $F$ with indices $\neq e$ (here $K_e$ is as in
Lemma \ref{lemma:Ken}), and
similarly $\Ltilde_e = K_e\cap\Ftilde$.
Now by Lemma \ref{lemma:Ken}, $K_e\cap F_e = K_e\cap\Ftilde_e=\Q$
for all $e$, so the splitting field
of $p_e(X)$ over $L_e$ has Galois group $\cS_7$, the symmetric
group on the seven roots of $p_e(X)$.  (It has this Galois group
over the larger ground field $K_e$ from the lemma, and
over $\Q$ itself, hence also over all intermediate fields, including $L_e$.)
This fact also shows that every automorphism $\sigma$ of $F$
fixes each $F_e$ setwise, since it must map $F_e$ into the intersection
of $F$ with the splitting field $P_e$ of $p_e(X)$ over $\Q$, and this
intersection is just $F_e$ itself.  (One says that $F_e$
is \emph{normal within $F$}.)

Now $\Gal{P_e}{\Q}$ is the symmetric group $\cS_7$ on the seven roots
of $p_e(X)$.  We know that $\Gal{P_e}{\Q(x_e+y_e)}$
contains exactly those permutations in $\cS_7$ that either interchange
$x_e$ with $y_e$ or fix both.  The $\Q$-conjugates of $(x_e+y_e)$ in $P_e$
are all sums of two distinct roots of $p_e$, and no other such sum
can be fixed by all these permutations, so $(x_e+y_e)$ has no
conjugates in $\Q(x_e+y_e)$.  If $(u_e+v_e)$ ever enters $F$,
then it is a conjugate of $(x_e+y_e)$, but then $x_e\in F$ as well,
so $\Gal{P_e}{F_e}$ then contains those permutations that fix the elements
$x_e$ and $y_e$ and the set $\{u_e,v_e\}$.  No other conjugate of
$(x_e+y_e)$ is fixed by all those permutations, and since there is
an element of $\Gal{P_e}{F_e}$ interchanging $u_e$ with $v_e$,
neither $u_e$ nor $v_e$ (nor any of the remaining three roots
of $p_e(X)$ in $\Qbar$) lies in $F_e$.

So,  if $\phi_e(x_e+y_e)\converges = x_e+y_e$,
then for every automorphism $\sigma$ of $F$,
we must have $\sigma(x_e)\in\{ x_e,y_e\}$, and hence
$\sigma(x_e+y_e)=x_e+y_e$, leaving $\sigma(u_e+v_e)=u_e+v_e$.
Thus $\sigma\res F_e$ has only two possibilities:  the identity,
and the map $\sigma_e\in\Aut{F_e}$ interchanging $x_e$ with $y_e$.  Note that
this $\sigma_e$ fixes $(x_e+y_e)$, and also fixes $(u_e+v_e)$.

On the other hand, if $\phi_e(x_e+y_e)$ diverges or converges
to any value $\neq x_e+y_e$, then $(x_e+y_e)$ has no conjugates
in $F$ except itself, and generates $F_e$.  So in this case $\sigma\res F_e$
must be the identity.

It follows that $B_F$ is computable.  First, of course, every pair $\la z,z\ra$
lies in $B_F$.  By the discussion above, $\la x_e+y_e,z\ra\in B_F$
if and only if $z=x_e+y_e$, and likewise $\la u_e+v_e,z\ra\in B_F$
if and only if $z=u_e+v_e$.
Next, for any root $z$ in $F$ of any $p_e$, we know
that $\phi_e(x_e+y_e)\converges=x_e+y_e$, since $F$ could
not contain such a root otherwise; hence
$\la z,z'\ra\in B_F$ if and only if $z'\in\{ x_e,y_e\}$.  Thus our decision
procedure for $B_F$ can compute
the orbit of every generator of $F_e$.  So, for an arbitrary pair
$\la z,z'\ra\in F^2$, it can express $z$ in terms of the generators
of finitely many $F_e$, compute the orbits of each of these
generators, and use this information to determine whether $\la z,z'\ra\in 
B_F$.
(We do use here the remark above that for any $e_0,\ldots,e_n$
with all $x_{e_i}\in F$, the choice of images of the $x_{e_i}$
may be made independently for all $i\leq n$.)

Because $F$ is a field, computability of $B_F$ implies
computability of $A_F$, as follows.  Given any tuples $\avec,\bvec\in F^n$,
the Effective Theorem of the Primitive Element allows us to identify
a single element $a\in F$ such that $\Q(a)=\Q(\avec)$.
Having found some $q\in\Q[X_1,\ldots,X_n]$ with $a=q(\avec)$,
we set $b=q(\bvec)$.  In turn, each $a_i=q_i(a)$
for some $q_i\in\Q[X]$.  Now $\la \avec,\bvec\ra\in A_F$
if and only if $\la a,b\ra\in B_F$ and, for each $i\leq n$, $b_i=q_i(b)$.
For full details, we refer the reader to \cite{MS11},
where the Effective Theorem of the Primitive Element
appears as Theorem 3.11.

However, for any $e$, we know that $\phi_e(x_e+y_e)\converges=x_e+y_e$
if and only if no isomorphism $\rho:F\to\Ftilde$ fixes $(x_e+y_e)$.  In particular, if
$\phi_e(x_e+y_e)$ diverges or converges to any other value,
then $(x_e+y_e)$ has no $\Q$-conjugate in either $F$ or $\Ftilde$;
whereas if $\phi_e(x_e+y_e)\converges=x_e+y_e$, then $x_e$ and $y_e$ lie in $F$,
but the only two roots of $p_e$ in $\Ftilde$ are $u_e$ and $v_e$,
whose sum is $\neq\phi_e(x_e+y_e)$.  In both these cases,
therefore, $\phi_e$ is not an isomorphism from $F$ onto $\Ftilde$.
Since this fact holds for all $e$, the field $F$ is not computably 
categorical.
\qed\end{proof}

\section{Computable Categoricity}
\label{sec:cc}

Relative computable categoricity immediately implies computable categoricity,
and for many theories $T$, the converse also holds of all computable models
of $T$.  This is the case for the theories of linear orders, Boolean algebras,
and trees (in the language of partial orders).  Indeed, in \cite{G75}
Goncharov showed
it to hold for all computable structures $\M$ for which the set of $\Sigma_2$
sentences in the elementary diagram of $\M$ is decidable.  On the other
hand, in \cite{K96}, Kudinov proved that decidability of the $\Sigma_1$
fragment of the elementary diagram does not suffice, by producing
examples of such \emph{$1$-decidable} structures that are computably
categorical but
not relatively computably categorical.  Such structures are widely
considered to be unusual, and among classes of structures for which
characterizations of computable categoricity are known (linear orders,
Boolean algebras, trees, etc.), computably categoricity always implies
relative computable categoricity.
Here we address this question for algebraic fields,
proving that computable categoricity of an algebraic field
can fail to relativize.

\begin{theorem}
\label{thm:notrcc}
There exists a computable, computably categorical algebraic field $F$
whose orbit relation $B_F$ is not $\Sigma^0_2$. It follows from
Theorem \ref{thm:rccorbit} that $F$ is not relatively computably categorical.
\end{theorem}
\begin{proof}
Every computable algebraic field has $\Pi^0_2$ orbit relation $B_F$,
by Lemma \ref{lemma:finiteblock}.  Therefore, it suffices to construct
such a field $F$ that is computably categorical, but for which
$B_F$ is not $\Delta^0_2$.

Our construction of $F$ takes place on a tree $T$, and satisfies
two types of requirements.  Every node $\beta$ at level $2e$ of $T$
is dedicated to satisfying requirement $\C_e$ for computable categoricity:
\begin{align*}
\C_e:~~&\text{If $\phi_e$ decides the atomic diagram of a field $C_e$
isomorphic to $F$,}\\
&\text{then there is a computable isomorphism $C_e\to F$.}
\end{align*}
Such a $\beta$ is also called a $\C_e$-node.  It has two outcomes $\cong$
and $\not\cong$, and hence two immediate successors, which we order:
$$\beta\hat{~}\la\cong\ra~~~~~~\prec~~~~~~\beta\hat{~}\la\not\cong\ra.$$
For the $\beta$ on the true path of $T$, the $\cong$ outcome
will hold infinitely often if and only if $C_e$ really is isomorphic to $F$,
in which case we must be sure to satisfy $\C_e$.
To guess at whether $C_e\cong F$, we use Corollary
\ref{cor:nofiniteblock}.  Each time we see both that a new larger subfield
of $C_e$ embeds into $F$, and that a new larger subfield of $F$
embeds into $C_e$, we make the outcome $\cong$ eligible.
The \emph{counter} $c_\beta$ is used for this purpose:
at stage $s$, it will be maximal such that the first $c_{\beta,s}$
elements of $C_{e,s}$ can be embedded together into $F_s$,
and the first $c_{\beta,s}$ elements of $F_s$ can
be embedded together into $C_{e,s}$.

Every node $\alpha$ at level $2e+1$ of the tree is dedicated to satisfying
the requirement $\R_e$ destroying limit-computability of $B_F$:
$$\R_e:~~\exists w \lim_t \phi_e(w,t)\neq B_F(w).$$
An $\R_e$-node $\alpha$ has only one outcome.  It acts so as
to satisfy its requirement by observing the behavior of its
function $\phi_e(\la x_\alpha,y_\alpha\ra,t)$ for a particular
pair $\la x_\alpha,y_\alpha\ra$ of elements of $F$, and ensures
that this pair lies in $B_F$ if and only if the limiting value of $\phi_e$ on the pair
says otherwise.

The construction somewhat resembles that of Theorem 4.1
of \cite{M09b}, which built a computable algebraic field
that is not $\bfz'$-categorical.  We will use here the same principal tool
given for that construction, and the same notation.
In particular, when $h(X)$ is a polynomial with coefficients in $\Q[\sqrtp]$,
we will write $h^-(X)$ to denote the image of this polynomial
under the nontrivial automorphism of $\Q[\sqrtp]$, with
$\sqrtp$ mapped to $-\sqrtp$. 

\begin{proposition}[Proposition 2.15 in Miller \cite{M09b}]
\label{prop:hpolys}
For any fixed prime $p$, let $F$ be the field $\Q[\sqrtp]$.
Then for every odd prime number $d$, there exists a polynomial
$h(X)\in F[X]$ of degree $d$ with the following properties.
\begin{itemize}
\item
$h$ and $h^-$ are both irreducible in the polynomial ring $F[X]$.
\item
The splitting field of $h$ over $F$ has Galois group isomorphic
to $\cS_d$, the symmetric group on the $d$ roots of $h$,
and the same holds for $h^-$.  (Since $\cS_d$ acts transitively on
the roots, this property implies the preceding one.)
\item
The splitting field of $h(X)$ over the splitting field of $h^-(X)$ also has
Galois group isomorphic to $\cS_d$ (and vice versa).
In particular, each of $h(X)$ and $h^-(X)$ is irreducible over the splitting field
of the other.
\end{itemize}
Moreover, uniformly in $p$, $d$, and any computable presentation of $F$,
it is computable whether an arbitrary $h(X)\in F[X]$ satisfies
these properties.
\end{proposition}

For an $\R_e$-node $\alpha$, $x_\alpha$ and $y_\alpha$
will be the two square roots of a specific prime number $p_\alpha$.
At each stage, $\alpha$ will use two polynomials $h_{b_\alpha-1,\alpha}$
and $h_{b_\alpha,\alpha}$, provided by the Proposition for this $p_\alpha$,
of distinct prime degrees, with $b_\alpha$ keeping count of these
polynomials.  Both $h_{b_\alpha-1,\alpha}$ and $h_{b_\alpha,\alpha}$
will have roots in $F$, but neither $h_{b_\alpha-1,\alpha}^-$ and $h_{b_\alpha,\alpha}^-$
will have a root.  Thus, $x_\alpha$ is tagged in two different ways
to distinguish it from $y_\alpha$.  Whenever
$\phi_e(\la x_\alpha,y_\alpha\ra,t)$ equals $0$ for a new larger $t$,
$\alpha$ takes a step toward making $\la x_\alpha,y_\alpha\ra$
lie in $B_F$:  it adjoins a root of $h_{b_\alpha-1,\alpha}^-$ to $F$, so that $y_\alpha$
now has this tag, just like $x_\alpha$.  At the same time, though,
$\alpha$ chooses a new $h_{b_\alpha+1,\alpha}$, giving $x_\alpha$
a new tag which $y_\alpha$ lacks, and increments $b_\alpha$ 
so as to keep track of the current tags.  Therefore, the only way
$\la x_\alpha,y_\alpha\ra$ can end up in $B_F$ is if this step
is repeated infinitely often, in which case $\lim_t\phi_e(\la x_\alpha,y_\alpha\ra,t)\neq 1$.
If the limit is $1$, then this step is repeated only finitely often,
and $x_\alpha$ is tagged in some way in which $y_\alpha$ never is,
so that $\la x_\alpha,y_\alpha\ra\notin B_F$.  The tagging of $x_\alpha$
by two separate polynomials ensures that higher-priority
$\C$-nodes can always distinguish $x_\alpha$ from $y_\alpha$,
so that they can build their computable isomorphisms, and that
once they have built them, their guesses in their fields $C_e$
will remain correct about which node corresponds to $x_\alpha$
and which to $y_\alpha$:  at least one tag will always be present,
keeping their computation correct, even while the other tag
is removed and redefined.

We start with $z_0$ as the multiplicative identity, so that $F_0$ is a copy of $\Q$.
All nodes are initialized
at stage $0$, which means that all counters $c_{\beta,0}$ for $\C$-nodes
$\beta$ and all primes $p_\alpha$, counters $b_\alpha$, and
potential field elements $x_{\alpha}$ and $y_{\alpha}$
for $\R$-nodes $\alpha$ are undefined.

To arrange our stages,
we fix the bijection $\omega\times\omega\to\omega$ according
to the listing of $\omega\times\omega$ as follows:
$$ \la 0,0\ra, \la 0,1\ra, \la 1,1\ra, \la 0,2\ra, \la 1,2\ra, \la 2,2\ra, \la 0,3\ra, \ldots.$$
At all stages $s+1$ with $s=\la 0,n\ra$, the root node is eligible.
At all stages $s+1$ with $s=\la m,n\ra$ and $m\neq 0$, the node
eligible at the preceding stage will have designated one of its
immediate successors (of length $m$) to be eligible at this stage.
Each eligible $\C$-node chooses a successor to be eligible,
while for an eligible $\R$-node, its unique successor is eligible.

At each stage $s+1$ where a $\C_e$-node $\beta$ is eligible, let $s'$
be the most recent stage (if any) at which $\beta$ was eligible.  If there 
has been 
no such stage $s'$ yet, or if $\beta$ has been initialized since
that stage, then we set $c_{\beta,s+1}=1$ and make
$\beta\hat{~}\la\not\cong\ra$ eligible at stage $s+1$.  Otherwise,
we consider the function $\phi_{e,s}$ as the (partial) characteristic function
for the atomic diagram of a structure $C_e$ in the language of fields.
If $\phi_{e,s}$ enumerates any facts inconsistent with the field axioms,
or if the initial segment of $\beta$ on which $\phi_{e,s}$ converges
is the same as that on which $\phi_{e,s'}$ converged,
then $\beta\hat{~}\la\not\cong\ra$ is eligible at substage $t+1$.
Otherwise, we consider the fragment $C_{e,s}$ of a field as described
by $\phi_{e,s}$.  Lemma \ref{lemma:frag} ensures that we
can check whether both of the following hold (as defined in Definition \ref{defn:frag}).
\begin{itemize}
\item
If $D$ is the finite subset $\{ 0,1,\ldots,c_{\beta,s'}\}$ of
the field fragment $\{0,1,\ldots,s\}$ of $F_s$, then $D$ embeds into $C_{e,s}$.
\item
The following finite subset $D'$ of the field fragment $C_{e,s}$ has an
embedding into the field fragment $\{ 0,\ldots,s\}$ of $F_s$.
The set $D'$ contains the elements $0,1,\ldots,c_{\beta,s'}$, and also
contains all elements $x\in C_{e,s}$ for which there exists
an $\alpha$ with $\beta\hat{~}\la\cong\ra\subseteq\alpha$ and $i\in\omega$
and a currently defined polynomial of the form $(X^2-p_{\alpha,s+1}$) or 
one of the forms
$h_{i,\alpha,s}$ or $h^-_{i,\alpha,s}$, such that when the coefficients
of this polynomial are mapped into $C_{e,s}$, $x$ is a root in $C_{e,s}$
of the resulting polynomial there.
(Notice that $(X^2-p_{\alpha,s+1})$ has coefficients in $\Q$,
hence has a unique image in $C_{e,s}$, or no image at all if the
necessary elements of $\Q$ have not appeared in $C_e$
by stage $s$.  The other two polynomials,
$h_{i,\alpha,s}$ and $h^-_{i,\alpha,s}$, have coefficients in
$\Q(p_{\alpha,s+1})$, but interchanging the two square roots
of $p_{\alpha,s+1}$ only interchanges these two
polynomials with each other.  So the description above
names at most three polynomials in $C_{e,s}[X]$ for each $\alpha$,
without ambiguity.)
\end{itemize}

If so, then we set $c_{\beta,s+1}=c_{\beta,s'}+1$,
initialize all nodes to the right of $\beta\hat{~}\la\cong\ra$
(including $\beta\hat{~}\la\not\cong\ra$ and all its successors),
and make $\beta\hat{~}\la\cong\ra$ eligible at stage $s+2$.
If either of these fragments fails to embed, then we leave
$c_{\beta,s+1}=c_{\beta,s'}$ and make $\beta\hat{~}\la\not\cong\ra$ eligible.

(The point of this stage is that, if $\beta$ is eligible
at infinitely many stages, then $\beta\hat{~}\la\cong\ra$
will be eligible at infinitely many stages if and only if $\phi_e$ is total and is the
characteristic function of the atomic diagram of a field $C_e$ such that
every finitely generated subfield of $C_e$ embeds into $F$ and vice versa.
By Lemma \ref{lemma:finiteblock}, this latter condition is equivalent to $C_e\cong F$.)

At each stage $s+1$ at which an $\R_e$-node $\alpha$ is eligible,
we let $s'$ be the most recent stage at which $\alpha$ was eligible
and executed either Step 1 or Step 3 (below).  In Step 1, $\alpha$
chooses its prime $p_\alpha$, adjoins its square roots $x_\alpha$
and $y_\alpha$, and sets up the first two tags on $x_\alpha$.
In Step 2, it waits for all $\C$-nodes $\beta$
with $\beta\hat{~}\la\cong\ra\subseteq\alpha$ to complete
their computations on the current tags, since $\alpha$ cannot add
any new tags until this has been done.  Finally, in Step 3,
$\alpha$ gets to check whether $\lim_t\phi_e(\la x_\alpha,y_\alpha\ra,t)$
has taken any further steps towards equalling $0$;
if so, then it executes the appropriate action with its tags,
while if not, it does nothing.
Here are the full descriptions of the steps.
\begin{enumerate}
\item
If this is the first stage at which $\alpha$ has been eligible
(so the stage $s'$ does not exist), or if $\alpha$ has been initialized
since stage $s'$, then we let $p_{\alpha,s+1}$ be the least
odd prime number not yet chosen as $p_{\gamma,t}$ at any previous stage
$t$ for any node $\gamma$.  We adjoin to $F_s$ a new element $x_{\alpha,s+1}$
satisfying $x_{\alpha,s+1}^2 = p_{\alpha,s+1}$.  Of course, this 
action also
adjoins a second square root of $p_{\alpha,s+1}$, and
we name this second element $y_{\alpha,s+1}$.
Then we choose the least two odd prime numbers not yet used in the construction,
and, for each of these two prime numbers $d$, search for a
polynomial as described in Proposition \ref{prop:hpolys}
of this degree $d$ for the prime $p_{\alpha,s+1}$.
Let $h_{0,\alpha,s+1}(X)$ and $h_{1,\alpha,s+1}(X)$ be the two polynomials
we find, and set $b_{\alpha,s+1}=1$.  We adjoin to $F_s(x_{\alpha,s+1})$
one root of $h_{0,\alpha,s+1}$ and one root of $h_{1,\alpha,s+1}$,
to form $F_{s+1}$.  This action completes this stage.

\item
Otherwise, $p_{\alpha,s'}$ and $b_{\alpha,s'}$ are already defined,
as are $x_{\alpha,s'}$ and $y_{\alpha,s'}$ and polynomials
$h_{0,\alpha,s'}(X),\ldots,h_{b_{\alpha,s'},\alpha,s'}(X)$.
We check whether, for every $i\leq e$, either the $\C_i$-node $\beta\subset\alpha$
has $\beta\hat{~}\la\not\cong\ra\subseteq\alpha$ or else the field
fragment $C_{i,s}$ contains a root of the minimal polynomial
of the least primitive generator of $F_{s'}$ over $\Q$.  If this is not the case,
then we do nothing at this stage.
If it is the case, then we execute Step 3 below.

\item
Let $u\geq 0$ be maximal with the property that all
$\phi_{e,s}(\la x_{\alpha,s'},y_{\alpha,s'}\ra,t)$ with $t<u$ converge,
and let $u'$ be the corresponding maximum with $s'$
in place of $s$.  If none of
$\phi_{e,s}(\la x_{\alpha,s'},y_{\alpha,s'}\ra,u'),\ldots,
\phi_{e,s}(\la x_{\alpha,s'},y_{\alpha,s'}\ra,u-1)$ equals $0$,
or if $u'=u$, then we do nothing.  Otherwise
we set $b_{\alpha,s+1}=b_{\alpha,s'}+1$,
choose the least odd prime $d$ not yet seen in the construction,
and find a polynomial $h_{b_{\alpha,s+1},\alpha,s+1}(X)$
satisfying Proposition \ref{prop:hpolys} for this $d$ and
for $p_{\alpha,s'}$.  We adjoin to $F_s$ one root of this
$h_{b_{\alpha,s+1},\alpha,s+1}$, and also one root
of $h^-_{b_{\alpha,s+1}-2,\alpha,s'}$.
(In Lemma \ref{lemma:irred} below, we show that
each of these polynomials is irreducible over $F_s$,
and indeed over the root of the other, so that this
suffices to define $F_{s+1}$.)  Of course, $b_{\alpha,s+1}-2=b_{\alpha,s}-1$,
so this has added a tag for $y_{\alpha,s+1}$ for that old $h$-polynomial,
as well as a tag for $x_\alpha$ for the new $h$-polynomial
$h_{b_{\alpha,s+1},\alpha,s+1}$.
\end{enumerate}
Whichever step was executed, we then end this stage,
with the unique successor of $\alpha$ eligible at the next stage.

Since this process was effective and every $F_{s+1}$ was an algebraic extension
(proper or not) of the preceding $F_s$, we have constructed a computable
algebraic field $F=\bigcup_s F_s$ of characteristic $0$.  We claim
that every requirement $\C_e$ and $\R_e$ is satisfied by our construction,
so that $F$ is the field needed to establish the theorem.

The $\C$-nodes $\beta$ always make one of their two successors eligible,
and $\R$-nodes $\alpha$ always make their unique successor eligible.
Therefore, the set containing the leftmost node at each level that
is eligible infinitely often forms a path through the tree, called
the \emph{true path} $P$.  Each requirement corresponds to a unique
node on $P$, which will be the node causing that requirement to be satisfied.
If $\beta\in P$ is a $\C$-node, then $\beta\hat{~}\la\not\cong\ra\in P$
if and only if $\lim_s c_{\beta,s}$ exists and is finite;
if $\beta\hat{~}\la\cong\ra\in P$, then $c_{\beta,s}\to\infty$.

We start with an analysis of the strategy of the $\R_e$-node
$\alpha\in P$, for any fixed $e$, starting with the last stage at which
this $\alpha$ is initialized.  Let $s_0$ be the first subsequent stage
at which it is eligible.  Then whenever $\alpha$ is eligible after $s_0$,
it has two particular $h_\alpha$ polynomials that have
roots, but such that the corresponding $h_\alpha^-$
polynomials do not have roots.  We say that $x_\alpha=x_{\alpha,s_0}$ is
\emph{tagged} by these polynomials, while $y_{\alpha,s_0}$ is
not (yet) tagged by them.  Suppose that $\alpha$ has the
correct guess about which $\C$-nodes preceding
it correspond to fields isomorphic to $F$,
and that $\phi_{e}(\la x_\alpha,y_\alpha\ra,t)\converges$
for all $t$ (since otherwise $\R_e$ is trivially satisfied).
If there are infinitely many
$t$ for which $\phi_{e}(\la x_\alpha,y_\alpha\ra,t) = 0$,
then the limit of $\phi_e$ on $\la x_\alpha,y_\alpha\ra$
can only equal $0$ or fail to exist,
yet both $x_\alpha$ and $y_\alpha$ do end
up each tagged by all the polynomials $h_{n,\alpha}$,
and thus will lie in the same orbit.  Conversely,
if there are only finitely many such $t$,
then $\lim_t\phi_{e}(\la x_\alpha,y_\alpha\ra,t)\neq 0$
(and the limit may not exist at all), but in this case
only finitely many polynomials $h_{n,\alpha,s}$ were ever
defined, and the last two still tag $x_\alpha$ without
tagging $y_\alpha$.  So in this case $\la x_\alpha,y_\alpha\ra\notin B_F$.
Thus in both cases, $\R_e$ will be satisfied.
Of course, to complete this argument, we must show that
the tags really do work the way we claimed here,
and in particular that no extraneous tags were introduced by
the actions of other nodes.

So consider the elements adjoined to $F_s$ at a specific
stage $s+1$, with an $\R$-node $\alpha$ eligible at this stage.
If $\alpha$ is in Step 1, then it first adjoins a square root $x_{\alpha,s+1}$
of its prime $p_{\alpha,s+1}$.  This extension has degree $2$,
since we chose a $p_{\alpha,s+1}$ that does not already have a square root 
in $F$.
Then the stage adjoins roots of $h_{0,\alpha,s+1}$ and $h_{1,\alpha,s+1}$,
which were chosen to be irreducible over $\Q(x_\alpha)$ with degrees
that are both new prime numbers, and indeed are irreducible
over $F_s$ as well, by Lemma \ref{lemma:irred} below.
Thus the root of $h_{i,\alpha,s+1}$
generates an extension of degree $\deg{h_{i,\alpha,s+1}}$, for each $i$,
and since these degrees are prime to each other, they generate
linearly disjoint field extensions of $F_s(x_{\alpha,s+1})$ (that is, field extensions
whose intersection equals just $F_s(x_{\alpha,s+1})$).  So the degree
$[F_{s+1}:F_s]$ is the product of these two primes and $2$.

Now suppose $\alpha$ is in Step 3 at stage $s+1$,
and adjoins to $F_s$ one root of $h_{b_{\alpha,s+1},\alpha,s+1}$ and
one root $r^-$ of $h^-_{\alpha,b_{\alpha,s+1}-2,s+1}$.  The former is irreducible
over $F_s$ and has a new large prime as its degree, and its root
thus generates an extension of that degree.  We also claim that
the root $r^-$ of $h^-_{\alpha,b_{\alpha,s+1}-2,s+1}$ generates
a further extension of degree $\deg{h^-_{\alpha,b_{\alpha,s+1}-2,s+1}}$.
The following lemma justifies these claims.
\begin{lemma}
\label{lemma:irred}
The following holds for every stage $s$.
First, for the $\alpha$ and $i$ (if any) such that 
$h_{i,\alpha,s+1}(X)$ is first defined at stage $s+1$,
the polynomial $h_{i,\alpha,s+1}(X)$ is irreducible in 
$F_s(x_{\alpha,s+1})[X]$.
Second, for any $\alpha$ and $i$ such that
$h_{i,\alpha,s}(X)$ has a root in $F_s$
but $h^-_{i,\alpha,s}(X)$ does not,
$h^-_{i,\alpha,s}(X)$ is irreducible in $F_s[X]$.

Indeed, at a stage at which $\alpha$ enumerates
two roots (apart from $x_{\alpha,s+1}$ itself) into $F_{s+1}$,
using Step 1 or Step 3, the minimal polynomial of each of these roots
(either $h_{i,\alpha,s+1}$ or $h^-_{i,\alpha,s+1}$)
remains irreducible over the extension of $F_s$
by the other root.
\end{lemma}
\begin{proof}
We prove the first two statements simultaneously,
by induction on $s$.  First suppose that $h_{i,\alpha,s+1}(X)$
is defined at stage $s+1$, by $\alpha$ in Step 1 or Step 3,
with a root $r$ adjoined to $F_s$.  By Proposition \ref{prop:hpolys},
$[\Q(x_{\alpha,s+1},r):\Q(x_{\alpha,s+1})]=d$, the degree of $h_{i,\alpha,s+1}(X)$,
and therefore $d$ divides $[F_s(x_{\alpha,s+1},r):\Q(x_{\alpha,s+1})]$.
However, the prime degree $d$ was never used for any $h$-polynomials
except $h_{i,\alpha,s+1}$.
Now we use our inductive hypothesis on previous stages,
noting that since the elements adjoined by $\R$-nodes
$\gamma$ at previous stages $t$ were roots of irreducible
polynomials $h_{j,\gamma,t}$ or $h^-_{j,\gamma,t}$
(or square roots of primes),
those adjoinments created extensions of prime degrees
distinct from $d$.  Therefore, $d$ must divide
$[F_s(x_{\alpha,s+1},r):F_s]$, and since $x_{\alpha,s+1}$ has degree
either $1$ or $2$ over $F_s$, we have that $d$ divides
$[F_s(x_{\alpha,s+1},r):F_s(x_{\alpha,s+1})]$.  On the other hand,
$r$ is a root of $h_{i,\alpha,s+1}(X)$, which has degree $d$,
and so $d=[F_s(x_{\alpha,s+1},r):F_s(x_{\alpha,s+1})]$,
forcing $h_{i,\alpha,s+1}(X)$ to be the minimal polynomial
of $r$ over $F_s(x_{\alpha,s+1})$.  Therefore $h_{i,\alpha,s+1}(X)$
is irreducible in $F_s(x_{\alpha,s+1})[X]$.

Next, suppose that $h_{i,\alpha,s+1}$ has a root in $F_{s+1}$ but
$h^-_{i,\alpha,s+1}$ does not.  Let $d=\deg{h^-_{i,\alpha,s+1}}$, and fix a root
$r\in F_{s+1}$ of $h_{i,\alpha,s+1}$.  Set $E=F_{s+1}(r^-)$, where $r^-$
is a root of $h^-_{i,\alpha,s+1}$.  (If $h^-_{i,\alpha,s+1}$ is reducible, then
$r^-$ may be a root of any of its irreducible factors in $F_{s+1}[X]$,
and the argument below will apply.)  $E$ thus contains roots of
both $h_{i,\alpha,s+1}$ and $h^-_{i,\alpha,s+1}$, and by Proposition
\ref{prop:hpolys}, $d^2$ must divide $[E:\Q(x_{\alpha,s+1})]$,
since $d$ divides both $[\Q(r,x_{\alpha,s+1}):\Q(x_{\alpha,s+1})]$
and $[\Q(r^-,r,x_{\alpha,s+1}):\Q(r,x_{\alpha,s+1})]$.
However, the prime degree $d$ was never used for any $h$-polynomials
except $h_{i,\alpha,s+1}$ and $h^-_{i,\alpha,s+1}$.
Using the inductive hypothesis once again,
we see that among all $[F_{t+1}:F_t]$ with
$t<s$, the only one divisible by $d$ is
the one for the stage $t+1$ with $r\in F_{t+1}-F_t$;
moreover, for this $t$, we have that $[F_{t+1}:F_t]$ is divisible by $d$
but not by $d^2$.
It follows that $d$ must divide $[E:F_s]$, and therefore
the minimal polynomial of $r^-$ over $F_s$ has
degree divisible by $d$.  But $r^-$ is a root
of $h^-_{i,\alpha,s+1}$, which itself has degree $d$,
and so this is the minimal polynomial of $r^-$ over $F_s$.
Thus $h^-_{i,\alpha,s+1}$ is irreducible over $F_s$.
This conclusion completes the induction.

Finally, considering the two roots enumerated into $F_{s+1}$
by $\alpha$, we note that their minimal polynomials over $F_s$
have distinct prime degrees.  Therefore, the field extensions
generated by each are linearly disjoint:  their intersection is $F_s$.
It follows that neither extension can cause the other minimal
polynomial to factor (see e.g.\ \cite[Lemma 2.12]{M09a};
or just extend the argument from the induction above).
\qed\end{proof}

\begin{corollary}
\label{cor:noroot}
Fix any $i$, $\alpha$, and $t$.  Then
$F$ has a root of $h^-_{i,\alpha,t}(X)$ if and only if the node $\alpha$
itself adjoins such a root by entering Step 3 at some stage
$s+1$ with $b_{\alpha,s+1}=i+2$ and with $h_{i,\alpha,t}=h_{i,\alpha,s+1}$.
\end{corollary}
\begin{proof}
By Lemma \ref{lemma:irred}, at any stage $s+1$, nodes $\gamma\neq\alpha$
enumerate only roots of polynomials that are irreducible
over $F_s$ and have prime degrees distinct from the degree
of $h_{i,\alpha,t}$.  The same holds for the node $\alpha$
itself at stages $s+1$ such that $\alpha$ is initialized between
that stage and stage $t$, or such that $\alpha$ adjoins
roots of polynomials $h_{j,\alpha,s+1}$ or $h^-_{j,\alpha,s+1}$
with $j\neq i$.  Finally, when $\alpha$ adjoins a root of
$h_{i,\alpha,t}$, Proposition \ref{prop:hpolys} shows that
no root of $h^-_{i,\alpha,t}$ can result.
\qed\end{proof}

To see that requirement $\R_e$ is satisfied, let $\alpha$ be the $\R_e$-node
on the true path $P$, i.e.\ the leftmost node at level $2e+1$ that
is eligible at infinitely many stages.  For simplicity, write $x$ and $y$ for
$\lim_s x_{\alpha,s}$ and $\lim_s y_{\alpha,s}$, write $h_i$ for $\lim_s h_{i,\alpha,s}$, etc.
Now for any $\C_e$-node $\beta\subseteq\alpha$,
either $\beta\hat{~}\la\not\cong\ra\subseteq\alpha$,
or else $\beta\hat{~}\la\cong\ra$ was eligible in between
every pair of stages at which $\alpha$ was eligible.
In the latter case, $\beta\hat{~}\la\cong\ra$ was eligible
infinitely often, and so the field $C_e$ must be isomorphic
to $F$.  Therefore, $\alpha$ cannot simply execute Step 2 at cofinitely many stages;
it executes Step 1 at the first stage at which it is eligible,
and enters Step 3 at infinitely many stages after that.

Suppose first that $\lim_t\phi_e(\la x,y\ra,t)$ exists and equals $1$.
Then at all but finitely many of the stages at which $\alpha$
goes through Step 3, it does nothing, and so we have a finite limit
$b=\lim_s b_{\alpha,s}$.  Set $d=\deg{h_b}$.
Now $h_b(X)$ has a root $r$ in $F$,
while $h^-_b(X)$ does not, by Lemma \ref{lemma:irred}.
Therefore $\la x,y\ra\notin B_F$, since any automorphism
of $F$ mapping $x$ to $y$ would have to map $r$
to some root of $h^-_b$.  So in this case $\R_e$ is satisfied.

On the other hand, if $\lim_t\phi_e(\la x,y\ra,t)=0$, then
$\alpha$ executed Step 3 infinitely many times,
and so $b_{\alpha,s}\to\infty$ as $s\to\infty$.
Thus, for every $b$, both $h_b$ and $h^-_b$ have
a root in $F$.  We claim that every subfield $F_s$ containing $x$
has an embedding into $F$ mapping $x$ to $y$.  Once this claim
is established, Lemma \ref{lemma:finiteblock}
will show that there exists a field embedding of $F$ into itself
mapping $x$ to $y$, which will prove that $\la x,y\ra\in B_F$, since by
Lemma \ref{lemma:endo} this embedding must be an automorphism.
The key to proving this claim is the following standard fact from field theory.
\begin{lemma}
\label{lemma:embed}
Let $K\subseteq L$ be a field extension generated by a single
$x\in L$  that is algebraic over $K$, and let $f:K\to E$
be a field embedding.  Fix the minimal polynomial $h\in K[X]$
of $x$ over $K$, and let $h^f$ be its image in $E[X]$
under the map $f$ on its coefficients.  Then $f$ extends
to an embedding of $L$ into $E$ if and only if $E$ contains a root of $h^f(X)$.
\end{lemma}

To prove the claim that every subfield $F_s$ containing $x$
has an embedding into $F$ mapping $x$ to $y$,
we show how to extend such embeddings from $F_s$ to $F_{s+1}$.
Assume by induction that $f$ is an embedding of $F_s$
into $F$ with $f(x)=y$ (and hence $f(y)=x$), with $f(x_{\gamma,t})=x_{\gamma,t}$
for all $\gamma\neq\alpha$ (and also with $\gamma=\alpha$ for
stages $t$ before the last initialization of $\alpha$), with $f(r)=r$ for 
every generator
$r$ adjoined by any $\R$-node $\gamma\neq\alpha$ or by $\alpha$
before its last initialization, and such that,
for every root $r\in F_s$ of any $h_{i,\alpha,s}(X)$, we have that
$f(r)$ is a root
of $h^-_{i,\alpha,s}(X)$, and for every root $r^-\in F_s$ of any
$h^-_{i,\alpha,s}(X)$, we have that $f(r^-)$ is a root of $h_{i,\alpha,s}(X)$.
(Of course, this assumption applies only to elements of $F_s$,
which might not include $x$ or others of the above.)
We may assume that some $\R$-node $\gamma$ is eligible at stage
$s+1$, since otherwise $F_{s+1}=F_s$.

Consider first the case where either $\gamma\neq\alpha$ or $\alpha$
is initialized after stage $s+1$.  Now $\gamma$
may adjoin $x_{\gamma,s+1}$ and a root of each of $h_{0,\gamma,s+1}$
and $h_{1,\gamma,s+1}$ to $F_s$.  If so, then $x_{\gamma,s+1}^2=p_{\gamma,s+1}$,
and each of these roots satisfies an irreducible polynomial
over $\Q(x_{\gamma,s+1})$.  So we set $f(x_{\gamma,s+1})=x_{\gamma,s+1}$,
which extends $f$ to $F_s(x_{\gamma,s+1})$, by Lemma \ref{lemma:embed},
and also set $f$ to be the identity on the roots of these $h$-polynomials
(which works for the same reason).  Alternatively, $\gamma$
may have been in Step 3 and have adjoined a root of
$h_{b_{\gamma,s+1},\gamma,s+1}(X)$ and a root of
$h^-_{b_{\gamma,s+1}-2,\gamma,s+1}(X)$.  By Lemma
\ref{lemma:irred}, both are irreducible over $F_s(x_{\gamma,s+1})$.
But by assumption, $f$ restricts to the identity on the
coefficients of both (which all lie in $\Q(x_{\gamma,s+1})$),
and so again we can extend $f$ to these roots
just by taking the identity map on them.

Now consider the case where $\gamma=\alpha$ and $\alpha$
is never initialized after stage $s+1$.  First,
if $\alpha$ adjoins $x_{\alpha,s+1}$ to $F_s$, then its negative
$y_{\alpha,s+1}$ also appears in $F_{s+1}$ and is conjugate
to $x_{\alpha,s+1}$ over $F_s$, so Lemma \ref{lemma:embed}
allows us to define
$f(x_{\alpha,s+1})=y_{\alpha,s+1}$.  In either Step 1 or Step 3,
whenever $\alpha$ defines a new polynomial
$h_{i,\alpha,s+1}(X)$ and adjoins a root $r$ of it,
we know (by our assumption that $b_{\alpha,s}\to\infty$)
that eventually $\alpha$ will also adjoin a root $r^-$
of $h^-_{i,\alpha,s+1}(X)$ to $F$, and so we define
$f(r)$ to be that $r^-$.  By Lemma \ref{lemma:embed},
this definition does extend $f$ to a field
embedding on $F_s(r)$.  Likewise, if $\alpha$ is in
Step 3 and adjoins some root $r^-$ of
$h^-_{b_{\alpha,s+1}-2,\alpha,s+1}$ to $F_s$, we know by Lemma
\ref{lemma:irred} that $h^-_{b_{\alpha,s+1}-2,\alpha,s+1}$
is the minimal polynomial of this $r^-$ over $F_s$,
and so it is safe to set $f(r^-)$ to equal $r$, since $r$
is a root of the image $h_{b_{\alpha,s+1}-2,\alpha,s+1}$
of $h^-_{b_{\alpha,s+1}-2,\alpha,s+1}$ under the map $f$
on its coefficients.  Thus in all cases we have extended
$f$ from $F_s$ to a field embedding of $F_{s+1}$ into $F$,
with $f(x)=y$ whenever $x\in F_s$.  This fact proves the claim,
and completes our argument that requirement $\R_e$ is satisfied.

Turning to the $\C$-requirements, for any $e$,
we let $\beta$ be the node of length $2e$ on $P$,
i.e.\ the leftmost node of that length
that is eligible at infinitely many stages.  Of course, $\beta$
works for the requirement $\C_e$.  Suppose first that $\beta\hat{~}\la\cong\ra$
is never eligible after some stage $s_0$.
Then $F$ contains only finitely many elements enumerated
by nodes $\alpha\supseteq\beta\hat{~}\la\cong\ra$,
and moreover $\lim_s c_{\beta,s}=c_{\beta,s_0}$ is finite.
But the elements $\{ 0,1,\ldots,c_{\beta,s_0}\}$
and those enumerated by these $\alpha$ together
generate a subfield of $F$ that must not embed into $C_e$,
or else the subfield of $C_e$ generated by $\{0,\ldots,c_{\beta,s_0}\}$
does not embed into $F$, since otherwise $\beta\hat{~}\la\cong\ra$
would have become eligible again.  So in this case $F\not\cong C_e$,
satisfying $\C_e$.

Therefore we may assume that $C_e\cong F$ and that
$\beta\hat{~}\la\cong\ra\in P$.  We construct a computable isomorphism
from $F$ onto $C_e$ as follows.  First, let $s_0$ be the last stage at which
any node to the left of $\beta$ is eligible.  We may start by assuming that we know
the restriction $g_{0}=f\res F_{s_0}$ of the given isomorphism $f:F\to C_e$,
since this knowledge requires only finitely much information, namely the
images of the finitely many generators of $F_{s_0}$.

Now let $s_1<s_2<\cdots$ be all stages $>s_0$
at which $\beta\hat{~}\la\cong\ra$ is eligible.
(We can compute this sequence, of course.)  We extend each $g_{n}$
to the finite field extension $F_{s_{n+1}}$ of $F_{s_n}$ in turn, as follows.
If an element $x$ was adjoined to $F$
by a node $\alpha$ to the right of $\beta\hat{~}\la\cong\ra$,
then $\alpha$ is initialized at stage $s_{n+1}$, so we simply check how many
elements that $\alpha$ enumerated into $F$ before stage $s_{n+1}$.
In particular, for any $s$ with $s_n<s<s_{n+1}$, let $t<s_{n+1}$
be the greatest stage before the next initialization of $\alpha$.
If $\alpha$ enumerated a root of some
polynomial $h_{i,\alpha,s}$, we check whether
it also enumerated a root of $h^-_{i,\alpha,s}$ by stage $t$ or not.
This will be the case for finitely many $i$, but eventually
we will reach an $i$ for which $F_{t}$ contains a root
of $h_{i,\alpha,t}$ but no root of $h^-_{i,\alpha,t}$.
(Indeed, the same holds for $i+1$ as well, since $\alpha$ always
keeps two tags on $x_\alpha$ which $y_\alpha$ does not yet have.)
Fixing this $i$, we find both square roots of $p_{\alpha,t}$ in $C_e$,
and find a root of $h_{i,\alpha,t}$ over one of those square roots;
we then map $x_{\alpha,s}$ to the conjugate with this root
(and the root of $h_{i,\alpha,s}$ to the root itself,
and likewise for $i+1$),
and map $y_{\alpha,s}$ to the other conjugate.
This mapping also then determines, for each $j$ such that the polynomials
$h_{j,\alpha,t}$ and $h^-_{j,\alpha,t}$ both have roots in $F_{s_{n+1}}$,
where these roots should be mapped.  

By our choice of $s_0$, the only other nodes that can
enumerate any element into $F$ between stages $s_n$ and $s_{n+1}$
are nodes $\alpha$ with $\beta\hat{~}\la\cong\ra\subseteq\alpha$.
So next we suppose that such an $\alpha$ enters Step 1
at stage $s+1$, with $s_n\leq s < s_{n+1}$, and adjoins
$x_{\alpha,s+1}$ and roots $r$ and $r'$ of the two polynomials
$h_{0,\alpha,s+1}$ and $h_{1,\alpha,s+1}$.
Recall that these polynomials both have coefficients in
the field $\Q(\sqrt{p_{\alpha,s+1}})$.
We wait for both square roots of $p_{\alpha,s+1}$
to appear in $C_e$, which must happen eventually,
since by assumption $F\cong C_e$.  Once they have appeared,
each one gives rise to an image in $C_e[X]$
of the polynomial $h_{0,\alpha,s+1}\in F_{s+1}[X]$,
since either square root of $p_{\alpha,s+1}$ can be used as
the square root in $h_{0,\alpha,s+1}$.  As soon as either of these
two polynomials in $C_e[X]$ acquires a root in $C_e$, we define
$g_{n+1}(r)$ to equal that root, and define $g_{n+1}(x_{\alpha,s+1})$
to equal the square root of $p_{\alpha,s+1}$ that gave rise to the
polynomial that has this root.  We also consider
the polynomial $h_{1,\alpha,s+1}(X)$ in $C_e[X]$ defined
using $g_{n+1}(x_{\alpha,s+1})$, and wait for this polynomial
to acquire a root in $C_e$, which then becomes the value of
$g_{n+1}(r')$.  All of these events must eventually happen, since $F\cong C_e$.
It remains to show that this definition of $g_{n+1}$ actually does
extend to an embedding of $F$ into $C_e$.

Notice first that at every subsequent stage $t$, the polynomial
$h_{b_{\alpha,t}-1,\alpha,t}$ has a root $r$ in $F_t$
and $h_{b_{\alpha,t},\alpha,t}$ has a root $r'$ there,
but neither $h^-_{b_{\alpha,t}-1,\alpha,t}$ nor $h^-_{b_{\alpha,t},\alpha,t}$
has any root in $F_t$.  This $b_{\alpha,t}$ stays fixed from one stage to the next
(starting with $b_{\alpha,s+1}=1$), except for stages at which
$\alpha$ enters Step 3.  At such stages, $C_e$ must contain images
$g_{n+1}(r)$ and $g_{n+1}(r')$, since we defined $g_{n+1}$ on $r$ and $r'$
as soon as we found those roots in $C_e$.  Also, no roots
of the $g_{n+1}$-images of $h^-_{b_{\alpha,t}-1,\alpha,s+1}$ and
$h^-_{b_{\alpha,t},\alpha,s+1}$ have appeared yet, because 
by the construction for the node $\beta$, such roots
would prevent $\beta\hat{~}\la\cong\ra$ from
becoming eligible (and so $F$ would never have enumerated
roots of $h^-_{b_{\alpha,t}-1,\alpha,s+1}$ and $h^-_{b_{\alpha,t},\alpha,s+1}$,
and thus $F$ would not have been isomorphic
to $C$).  In Step 3, $\alpha$ adjoins to $F_s$ a root
of $h^-_{b_{\alpha,t}-1,\alpha,t}$, but $h^-_{b_{\alpha,t},\alpha,t}$
still has no root, and a new polynomial $h_{b_{\alpha,t}+1,\alpha,t}$
is defined, with a root $r''$ in $F_{t+1}$ but such that
$h^-_{b_{\alpha,t}+1,\alpha,t}$ has no root there.
So the situation remains the same, except that one
of the two holding polynomials has been replaced by a new one.
Before $\alpha$ can be eligible again, $C_e$ must acquire
an image for $r''$, but cannot acquire any root for
the $g_{n+1}$-image of $h^-_{b_{\alpha,t},\alpha,t}$.

Having understood the above, we consider three cases.
\begin{enumerate}
\item
If $\alpha$ is initialized at some stage $t+1>s+1$, then
the $g_{n+1}$-image of $h_{b_{\alpha,t},\alpha,t}$ has a root
in $C_e$, but the $g_{n+1}$-image of $h^-_{b_{\alpha,t},\alpha,t}$
will never acquire one.  Therefore, our choice of $g_{n+1}(x_{\alpha,s+1})$
was correct.
\item
If $\alpha$ is never initialized after stage $s+1$ but is only eligible
at finitely many stages, let $t$ be the last stage at which it is eligible.
The exact same analysis applies here as for the case when $\alpha$ is re-initialized.
\item
Otherwise $\alpha$ is never again initialized, but
is eligible infinitely often.  In this case,
$\alpha$ must enter Step 3 infinitely many times
(since $F\cong C_e$ precludes it from staying in Step 2 forever),
and so $\la x_\alpha,y_\alpha\ra\in B_F$, as discussed above.
Therefore, either of the square roots of $p_\alpha$
in $C_e$ can be the image of $x_\alpha$ under an isomorphism.
So in this case either choice for $g_{n+1}(x_{\alpha,s+1})$
would have been correct.
\end{enumerate}

Thus our definition of $g_{n+1}$ was correct for every node
$\alpha$ going through Step 1.  It remains
to define $g_{n+1}$ on all elements adjoined
to $F$ at any stage $s+1$ between stages $s_n$ and $s_{n+1}$ by $\R$-nodes
$\alpha$ in Step 3 of the construction.  But this definition is simple,
because for such an $\alpha$, the value $x_{\alpha,s+1}$ must already
have been defined, and we have already defined $g_{n+1}(x_{\alpha,s+1})$.
Therefore, when $h_{b_{\alpha,s+1},\alpha,s+1}$ is given a root $r$ in $F$,
we know the image of $h_{b_{\alpha,s+1},\alpha,s+1}$ in $C_e[X]$ under
the map $g_{n+1}$ on its coefficients, and we wait for this
image to acquire a root in $C_e$, which then becomes
$g_{n+1}(r)$.  Likewise, we know the image of
$h_{b_{\alpha,s+1}-2,\alpha,s+1}(X)$ under $g_{n+1}$,
and so we may wait for it to acquire a root in $C_e$,
then define $g_{n+1}$ to map the root of $h_{b_{\alpha,s+1}-2,\alpha,s+1}$
in $F$ to this root in $C_e$.  Since $C_e\cong F$, and since
$g_{n+1}(x_{\alpha,s+1})$ is correctly defined, such roots must appear.

Thus we have extended $g_{n+1}$ to all elements adjoined by any $\R$-node
between stages $s_n$ and $s_{n+1}$, so we have defined our
computable embedding $g_{n+1}$ on all of $F_{s_{n+1}}$.
It is clear that this process can then continue to $F_{s_{n+2}}$
and beyond, so that $g=\bigcup_n g_n$ is a computable embedding of $F$
into $C_e$.  But since we know that $C_e\cong F$, Lemma
\ref{lemma:endo} shows that $g$ is a computable isomorphism,
and so $\C_e$$\Sigma^0_1$  is satisfied.

The satisfaction of the requirement $\R_e$ shows that
$\lim_s\phi_e(~\!\cdot~\!,s)$ is not the characteristic function
of $B_F$, and so all these requirements together prove that
$B_F$ is not $\Delta^0_2$, hence not $\Sigma^0_2$.
On the other hand, the satisfaction of the $\C$-requirements
shows that $F$ is computably categorical, since every
computable field isomorphic to $F$ has an atomic diagram decidable
by some $\phi_e$, meaning that it is the field $C_e$, which
was made computably isomorphic to $F$ by the requirement $\C_e$.
These conclusions complete the proof of Theorem \ref{thm:notrcc}.
\qed\end{proof}

\section{Complexity of Computable Categoricity}
\label{sec:complexity}

Ostensibly, computable categoricity is a $\Sigma^1_1$ property,
since its definition involves the existence of (classical)
isomorphisms, hence involves quantifying over functions
from $\omega$ to $\omega$.  However, for those classes of
structures for which an exact complexity is known,
it has always turned out to be far less complex than $\Sigma^1_1$.
For instance, a computable linear order $\L$ is computably categorical
if and only if $\L$ contains only finitely many pairs of adjacent points,
and this condition can be expressed as a
$\Sigma^0_3$ formula in the (computable) order relation
on $\L$.  Indeed, for arbitrary computable structures $\M$,
the statement ``$\M$ has a $\Sigma^0_1$ Scott family'' is $\Sigma^0_3$,
and so relative computable categoricity is always a
$\Sigma^0_3$ property.

For algebraic fields, the very fact of being isomorphic is nowhere
near $\Sigma^1_1$.  Corollary \ref{cor:nofiniteblock} shows
that for algebraic fields $E$ and $F$,
being isomorphic is $\Pi^0_2$, since for any finitely generated
subfield $F_0$ we can effectively find a primitive generator of $F_0$,
and then find the minimal polynomial of that generator over the prime subfield
of $F_0$, so that the embeddability of $F_0$ into $E$ reduces
to the existence in $E$ of a root of that minimal polynomial
(translated from the prime subfield of $F$ to that of $E$, of course).
Thus, algebraic fields $E$ and $F$ over the same prime subfield $Q$
are isomorphic if and only if
$$ (\forall p(X)\in Q[X])\left[(\exists x\in E~p(x)=0)\phantom{\frac12}\!\!\iff\phantom{\frac12}\!\!
(\exists y\in F~p(y)=0)\right].$$
If we write $C_e$ for the field (if any) whose atomic diagram
has characteristic function $\phi_e$, as in the proof of Theorem
\ref{thm:notrcc}, then we can discuss various complexities exactly.
\begin{proposition}
\label{prop:complexity}
All of the following sets are $\Pi^0_2$-complete.
\begin{itemize}
\item
$\Fld = \set{e}{C_e\text{~is a field}}$.
\item
$\textbf{AlgFld}~\!=\set{e}{C_e\text{~is an algebraic field}}$.
\item
$\set{\la e,i\ra}{C_e\text{~and $C_i$ are isomorphic algebraic fields}}$.
\item
$\set{i}{C_i\text{~is isomorphic to the field~}C_e}$, where $C_e$
is any fixed algebraic field.
\end{itemize}
\end{proposition}
\begin{proof}
$\Pi^0_2$ definitions of all these sets except $\Fld$
are readily produced, given Corollary \ref{cor:nofiniteblock}
and our discussion above.  Saying that $\phi_e$ is the
characteristic function of the atomic diagram of a field
requires saying that $\phi_e$ is total with range $\{ 0,1\}$
(a $\Pi^0_2$ property) and that the field axioms are satisfied by this diagram.
As usually stated, most of the field axioms are $\Pi^0_2$, but
the existence of an identity element for each operation appears
to be $\Sigma^0_2$, and the existence of inverses (stated below
for multiplication) appears to be $\Sigma^0_3$:
$$ \exists c\forall x\exists y (x+x=x~\text{or}~x\cdot y=y\cdot x=c).$$
This sentence can be reduced to a $\Pi^0_2$ statement simply by having
constant symbols for $0$ and $1$ in the signature, but it is worth
noting that even without such constant symbols, the field
axioms are still $\Pi^0_2$.
\begin{lemma}
\label{Pi2axioms}
A structure $\M$ in the signature with $+$ and $\cdot$ is a field
if and only if these two operations are both associative and commutative,
$\cdot$ distributes over $+$, and the following hold:
$$ \exists x\exists y (x\neq y)~~\&~~\forall x\forall y \exists u (x+u=y)~~\&~~
\forall x\forall y\exists u(x+x\neq x\implies x\cdot u=y).$$
Thus the field axioms can be expressed as a single first-order
$\forall\exists$ sentence.
\end{lemma}
\begin{proof}
The forwards implication is immediate, so assume that the given axioms hold.
Fix any single $x$, and apply the middle axiom to get a $u$ with $x+u=x$.
But now for any $y$, we have some $v$ with $x+v=y$ and
hence, given associativity and commutativity,
$$ y+u = (x+v)+u = v+(x+u) = v+x =y,$$
so that this $u$ is actually an additive identity element $0$.
The given axiom for addition then yields additive inverses.
But once we have these, we see that $x+x=x$ implies $x=0$,
so there must exist a $y$ with $y+y\neq y$ (lest $\M$
have only one element).  Then we repeat for multiplication
the same argument as for addition, using this $y$
to get the identity element.
\qed\end{proof}
The $\Pi^0_2$-completeness of the sets in Proposition \ref{prop:complexity}
is mostly an elementary exercise.  One easily shows that $\Inf=\{e:|W_e|=\infty\}$ (where $W_e$ is the
$e$th c.e.\ set) is
$1$-reducible to $\Fld$, for instance, just by fixing a
computable field $F$ and, on a given input $e$, building the
characteristic function
of the decision procedure for $F$ one element at a time, each time
we get further evidence that $e \in \Inf$ (i.e., each time a new
element enters $W_e$).  It is worth noting that
each of the other three sets
is $\Pi^0_2$-complete (under $1$-reductions) within the class $\Fld$.
(The relevant definition can be found in \cite[Defn.\ 1.2]{CHKM06}.)
For instance, there is a computable injective function $f$ such that
$\forall e (f(e)\in\Fld)$, but the field $C_{f(e)}$ is algebraic
if and only if $e\in\Inf$.  (Start building the field $\Q(X_0,X_1,\ldots)$
of infinite transcendence degree, one element at a time, and when
$e$ gets its $n$-th chip, turn $X_n$ into a rational number itself,
so large that it has not yet been ruled out by the finitely many
elements currently in $C_{f(e)}$.)
\qed\end{proof}

Since classical isomorphism is so easily expressed for algebraic
fields, the complexity of computable categoricity for the class
becomes much simpler than $\Sigma^1_1$.

\begin{proposition}
\label{prop:pi4}
For algebraic fields, the property of being computably categorical
is $\Pi^0_4$.
\end{proposition}
\begin{proof}
We simply write out the definition of computable categoricity
and apply Proposition \ref{prop:complexity}.  The computable algebraic
field $F=C_e$ is computably categorical if and only if:
$$ (\forall i)[(i\in\Fld~\&~ C_i\cong C_e)\implies \exists j
(\phi_j\text{~is an isomorphism}:C_i\to C_e)].$$
The statements $i\in\Fld$ and $C_i\cong C_e$ are both $\Pi^0_2$.
For $\phi_j$ to be an isomorphism, it must be total
(which is $\Pi^0_2$) and must preserve the field structure:
$$ \forall x\forall y [\phi_j(x+y)=\phi_j(x)+\phi_j(y)~~~~~\&~~~~~
\phi_j(x\cdot y)=\phi_j(x)\cdot\phi_j(y)],$$
which is $\Pi^0_1$ once we know that $\phi_j$ is total.
(For $\phi_j$ to have image $\omega$ is also $\Pi^0_2$,
but in fact is not needed here, by Corollary \ref{cor:endo}.)
\qed\end{proof}

Our main theorem for this section is the complementary property:
that for computable algebraic fields, computable categoricity is $\Pi^0_4$-hard,
and therefore $\Pi^0_4$-complete.  This theorem is substantially different
from previously known results about the complexity of
computable categoricity for specific classes of structures,
and thus serves to distinguish algebraic fields
from all those other classes.  In particular, all previously
known cases were $\Sigma^0_n$-complete for some $n$,
usually for $n=3$, so $\Pi^0_4$-completeness suggests
that something distinctly different is happening here.

\begin{theorem}
\label{thm:Pi04complete}
For computable algebraic fields, the property of being computably
categorical is $\Pi^0_4$-complete.
\end{theorem}
\begin{proof}
With Proposition \ref{prop:pi4} already proven, it remains
to show hardness.  Let $S$ be any $\Pi^0_4$-complete set,
such as the complement of $\emptyset^{(4)}$.
Since the set $\textbf{Inf}$ is $\Pi^0_2$-complete, we may express
$S$ by fixing some 1-1 total computable function $f$ for which:
$$ S =\set{n\in\omega}{\forall a\exists b (f(n,a,b)\in\textbf{Inf})}
=\set{n}{\forall a\exists b~|W_{f(n,a,b)}|=\infty}.$$
It will simplify our construction to assume that every
set $W_{f(n,a,b)}$ contains the element $0$,
and that at each single stage, at most one set
$W_{f(n,a,b)}$ receives a new element.

We will describe a 1-1 total computable function that maps
each $n\in\omega$ to the index for some computable algebraic field
$F$, which will be computably categorical if and only if $n\in S$.
The output of this function is the program that uses the following 
construction
(which is uniform in $n$) to build a computable field.  At the end of the construction,
we will demonstrate that the computable algebraic field $F$ that it built
is computably categorical if $n\in S$, but not otherwise.

The construction of $F$ is performed on a tree $T$, in a style reminiscent of that in
the proof of Theorem \ref{thm:notrcc}, adapted to incorporate the question
of whether $\forall a\exists b~f(n,a,b)\in\textbf{Inf}$.  As there, we let
$C_e$ denote the structure (in the language of fields) whose atomic
diagram is decided by the partial function $\phi_e$.  The tree $T$
for the construction will consist of two types of nodes.  We now
describe the basic modules used by each type to satisfy its requirement.

Every node $\beta$ at level $2e$ of $T$ 
is a \emph{categoricity node}, or \emph{$\C$-node}, dedicated
to satisfying requirement $\C_e$ for computable categoricity for $F$:
$$
\C_e:~C_e\cong F \implies\exists\text{~a computable isomorphism~}g_e:C_e\to F.
$$
Such a $\mathcal{C}_i$-node $\beta$ has two outcomes, 
$\cong$ and $\not\cong$, ordered with $\cong~\!\!\prec~\!\!\fin$.
The outcome $\cong$ denotes that the hypothesis of $\C_e$ turned
out to be true:  $C_e\cong F$.
In this case, the
$\beta$ on the true path at level $2e$ will produce the computable isomorphism
$g_\beta$ required, since no node above it or to its right will ever add
anything to $F$ that could cause problems for its isomorphism.
This process is much the same as that performed by the categoricity nodes
in the tree for Theorem \ref{thm:notrcc}.
The outcome $\fin$ denotes the negation of the outcome $\cong$,
in which case $\mathcal{C}_e$ holds automatically.

Every node $\alpha$ at level $2a+1$ of the tree is a \emph{non-categoricity node},
or \emph{$\R$-node}, trying to construct a computable field
$E_\alpha\cong F$ to satisfy the opposite requirement:
$$\R_a:~~[\forall b~f(n,a,b)\notin\Inf]\implies
[\forall b~\phi_b:E_\alpha\to F\text{~is not an isomorphism}].$$
The construction will build the computable fields $E_\alpha$
for every $\R$-node $\alpha$, all isomorphic to $F$.
An $\R_i$-node $\alpha$ has outcomes ordered in order type $\omega$:
$$ 0~\prec~1~\prec~2~\prec\cdots$$
If $\alpha$ lies on the true path, then for the least $b\in\omega$
(if any) such that $f(n,a,b)\in\Inf$, the node $\alpha\hat{~}\la b\ra$
will be the leftmost successor eligible infinitely often.
If there is no such $b$, then the hypothesis
of $\R_a$ is satisfied, and in fact the true path will end at $\alpha$;
in this case, the field $E_\alpha$ built by $\alpha$ will prove
that $F$ is not computably categorical.

The $\R_i$-node $\alpha$ runs the following basic module
simultaneously for all $b\in\omega$, although
whenever $W_{f(n,a,b)}$ receives a new element, $\alpha$
restarts its strategy for every $b'\geq b$.  For each $b$,
$\alpha$ starts by adjoining one witness element to $E_\alpha$
(with a corresponding witness adjoined to $F$)
and waits for $\phi_b$ to map the witness in $E_\alpha$
to the witness in $F$, which is its unique possible image there.
If $\phi_b$ does so, then $\alpha$ adds
a new element to $F$ to ``tag'' the witness there.
It waits until all categoricity nodes $\beta$
with $\beta\hat{~}\la\cong\ra\subseteq\alpha$ have mapped
the witness and its tag to an appropriate image, then adjoins a second
witness to $F$, conjugate to the original witness there,
and likewise adjoins a second witness to $E_\alpha$.
However, in $E_\alpha$, $\alpha$ tags the second witness
instead of the first.  Therefore, assuming no further
tags nor conjugates of the two witnesses ever appear in $F$,
$\phi_b$ cannot be an isomorphism, since it mapped the untagged
witness in $E_\alpha$ to the tagged witness in $F$.

All through this process (and forever after), $\alpha$
keeps watching to see if $W_{f(n,a,b)}$ receives any more elements.
If it ever does, then $\alpha$ terminates its procedure for $b$
and for all $b'>b$, makes $\alpha\hat{~}\la b\ra$
eligible and begins its entire process over again with a new
witness (which is the root of a completely new minimal polynomial).
Therefore, $\alpha$ precludes $\phi_b$ from being an isomorphism
only if $f(n,a,b)\notin\Inf$.  If every $f(n,a,b)\notin\Inf$, then all
of $\alpha$'s basic modules succeed, leaving $E_\alpha$
isomorphic to $F$ but not computably isomorphic to it.

At stage $0$, we begin with $F_0=\Q$ and also all fields $E_{\alpha,0}=\Q$.
All nodes are initialized, so that all values mentioned below for each node
are undefined at stage $0$.

The stages are ordered as in the construction in Theorem
\ref{thm:notrcc}, so that the root is eligible at every stage $\la 0,k\ra+1$,
and at each stage $\la l,k\ra+1$, some node at level $l$ is eligible
and (if $l<k$) chooses a node at level $(l+1)$ to be eligible
at the following stage $\la l+1,k\ra+1$.

At stage $s+1$, suppose that the $\mathcal{C}_e$-node
$\beta$ is eligible.  Let $s'$ be the greatest stage $\leq s$ at which
either $\beta$ was initialized or the node $\beta\hat{~}\la\cong\ra$
was eligible.  If the length of agreement between $F_s$ and $C_{e,s}$ (as defined
in the proof of Theorem \ref{thm:notrcc}) is no greater than
the domain of $g_{\beta,s'}$,
then we do nothing at this stage, and make $\beta\hat{~}\la\fin\ra$
eligible at the next stage.
If the length of agreement has increased, then $\beta\hat{~}\la\cong\ra$
will be eligible at the next stage.  At this stage, we define the
map $g_{\beta,s+1}$ to extend
the map $g_{\beta,s'}$ to the next element of the field $C_e$.
(By assumption, this must be a partial field embedding.)
This completes the stage.

At a stage $s+1$ at which an $\mathcal{R}_a$-node $\alpha$ is eligible,
we again let $s'$ be the greatest stage $\leq s$ at which
$\alpha$ either was initialized or was eligible.  Fix the least $b_0$ for which
$W_{f(n,a,b_0),s'}\neq W_{f(n,a,b_0),s+1}$.  (If there is no such $b_0$, then
find the least $t>s+1$ for which $(\exists b_0)W_{f(n,a,b_0),s'}\neq W_{f(n,a,b_0),t}$,
and choose that $b_0$.  Since all the sets
$W_{f(n,a,b)}$ are nonempty and only one can receive
an element at any given stage, we eventually find such a $b_0$.)
The node $\alpha\hat{~}\la b_0\ra$ will be eligible at the next stage.

If $\alpha$ was initialized at stage $s'$, then we simply set
both $E_{\alpha,s+1}$ and $F_{s+1}$ to equal $F_s$,
and end this substage.  If $\alpha$ was not initialized at stage $s'$,
then we execute the following instructions.

For each $b\geq b_0$, we initialize the $\alpha$-strategy for $b$,
by making $p_{\alpha,b,s+1}$ and all related roots, witnesses, and tags undefined.
First, however, for each $b\geq b_0$ for which $\alpha$ is
currently waiting to perform Step 3 (so $x_{\alpha,b,s'}\in F_s$,
but $E_{\alpha,s}$ does not yet contain any element $u_{\alpha,b,s'}$),
we adjoin $\xtilde_{\alpha,b,s'}$ to $E_{\alpha,s}$ (and then make
$\xtilde_{\alpha,b,s+1}$ undefined, along with all other roots and tags).
This ensures that $E_{\alpha,s+1}$ becomes isomorphic to $F_{s+1}$
once again (except possibly for certain tags for $\alpha$-strategies for values
$b'<b_0$; such tags might still lie in $F_s$ but have no images in $E_{\alpha,s+1}$).

For each $b<b_0$, we proceed according to the following steps.

\begin{enumerate}
\item
If no polynomial $p_{\alpha,b,s'}(X)$ is currently defined,
then we use Proposition \ref{prop:Ken} to choose a polynomial
$p_{\alpha,b,s+1}(X)\in\Q[X]$ of degree $7$, whose Galois group
(over the splitting field of the product of all $p$-polynomials used so far
in the construction, i.e.\ all $p_{\alpha',b',t}(X)$ with $t\leq s$)
is the symmetric group $\cS_7$ on its seven roots.
(Here we regard $\Q$ as a subfield of $F_s$, so that this polynomial
lies in $F_s[X]$.)  We define $x_{\alpha,b,s+1}$ and $y_{\alpha,b,s+1}$
to be two roots of $p_{\alpha,b,s+1}(X)$, but at this step we only
adjoin their sum $(x_{\alpha,b,s+1}+y_{\alpha,b,s+1})$ to $F_s$, forming
$F_{s+1}$ and leaving the roots themselves for possible later use.
Likewise, we adjoin the sum $(\xtilde_{\alpha,b,s+1}+\ytilde_{\alpha,b,s+1})$
of two roots of $\ptilde_{\alpha,b,s+1}(X)$ to every field $E_{\alpha',s}$,
for every $\R$-node $\alpha'$ including $E_{\alpha,s}$ itself.
So each such field $E_{\alpha',s+1}$ remains isomorphic to
$F_{s+1}$ (unless $E_{\alpha',s}\not\cong F_s$).
We define $q_{\alpha,b,s+1}(X)\in\Q[X]$ to be the minimal polynomial
of $(x_{\alpha,b,s+1}+y_{\alpha,b,s+1})$ over $\Q$.
Roots of $q_{\alpha,b,s+1}(X)$ will be called \emph{witnesses} being
used for $\alpha$ and $b$, in their respective fields $F$ and $E_{\alpha}$.
(In Step 3, a second witness may be adjoined to each of $F$ and $E_\alpha$.)
\item
If $x_{\alpha,b,s'}$ and $y_{\alpha,b,s'}$ are already defined
but $(x_{\alpha,b,s'}+y_{\alpha,b,s'})$
has not yet been tagged (as below), then we check whether
$\phi_{b,s}(x_{\alpha,b,s'}+y_{\alpha,b,s'})\converges = 
(\xtilde_{\alpha,b,s+1}+\ytilde_{\alpha,b,s+1})$.  If not,
then we do nothing at this stage.  If so, then we adjoin $x_{\alpha,b,s'}$
to $F_s$, calling it a \emph{tag} for the witness $(x_{\alpha,b,s'}+y_{\alpha,b,s'})$.
To preserve isomorphisms, we also adjoin $\xtilde_{\alpha,b,s'}$
to $E_{\alpha',s}$ for every $\R$-node $\alpha'$ except $\alpha$,
keeping $E_{\alpha',s+1}\cong F_{s+1}$ (unless $E_{\alpha',s}\not\cong F_s$).
Thus we leave $E_{\alpha,s+1}=E_{\alpha,s}\not\cong F_{s+1}$,
with no tag adjoined to $E_{\alpha,s}$.
\item
If $x_{\alpha,b,s'}\in F_s$ already, and $E_{\alpha,s'}$ contains
no corresponding tag, then we check whether,
for every $\C_e$-node $\beta$ with $\beta\hat{~}\la\cong\ra\subseteq\alpha$,
the domain of $g_{\beta,s+1}$ contains $x_{\alpha,b,s'}$
and the field fragment $C_{e,s+1}$
contains exactly one witness for $\alpha$ and $b$.
If not, then we do nothing.  If so, then we define
$u_{\alpha,b,s+1}$ and $v_{\alpha,b,s+1}$ to be new roots
of $p_{\alpha,b,s+1}(X)$, adjoin their sum $(u_{\alpha,b,s+1}+v_{\alpha,b,s+1})$
to $F_s$ as a new witness, and likewise
adjoin a new witness $(\utilde_{\alpha,b,s+1}+\vtilde_{\alpha,b,s+1})$,
the sum of two new roots of $\ptilde_{\alpha,b,s'}(X)$, to every $E_{\alpha',s}$
with $\alpha'\neq\alpha$.  To $E_{\alpha,s}$ we adjoin the two new roots
$\utilde_{\alpha,b,s+1}$ and $\vtilde_{\alpha,b,s+1}$ of $\ptilde_{\alpha,b,s'}(X)$;
this also adjoins their sum, of course, as a new witness,
and leaves $F_{s+1}\cong E_{\alpha,s+1}$,
but only via isomorphisms mapping $(x_{\alpha,b,s'}+y_{\alpha,b,s'})$
to $(\utilde_{\alpha,b,s+1}+\vtilde_{\alpha,b,s+1})$, since these are
the witnesses in their respective fields that now have tags.
This situation will be preserved forever (unless either $\alpha$ is initialized
or some $W_{f(n,a,b')}$ with $b'\leq b$ later receives a new element),
and so $\phi_b$ cannot be an isomorphism from $F$ onto $E_{\alpha}$.
\item
If none of the foregoing conditions applies, then $\alpha$ has satisfied
$\R_a$, and we do nothing at this stage.
\end{enumerate}
Having completed these steps for every $b<b_0$, we have finished
this stage.

When stage $s+1$ is completed, we initialize every
node to the right of the node eligible at that stage
(exactly as we did for $\alpha$-strategies for each $b\geq b_0$
in the construction for $\R$-nodes).  For a $\C$-node $\beta$, initialization
simply means that $g_{\beta,s+1}$ becomes the empty function.
For an $\R$-node $\alpha$, and for every $b\in\omega$,
we make all polynomials, roots, and tags associated with
$\alpha$ undefined at stage $s+1$, and we also make
$E_{\alpha,s+1}$ undefined.  This completes stage $s+1$.

It is clear that this construction builds a computable algebraic field $F$,
uniformly in $n$, and that this field is the extension of $\Q$ generated
by various witnesses and tags adjoined by assorted $\R$-nodes.
We claim that $F$ is computably categorical if and only if $n\in S$,
which is to say, if and only if for every $a$ there is some $b$ with $f(n,a,b)\in\Inf$.
As usual, the proof is based on the \emph{true path} $P$ through $T$,
i.e.\ the set of all nodes in $T$ that are eligible at infinitely
many stages, but initialized only finitely many times.

Suppose first that $n\in S$.  Now every $\C$-node $\beta$ makes
one of its two successors eligible whenever $\beta$ itself
is eligible.  Moreover, an $\R_a$-node $\alpha$ on $P$ will make
its successor $\alpha\hat{~}\la b\ra$ eligible infinitely often,
where $b$ is minimal such that $f(n,a,b)\in\Inf$, while for every
$b'<b$, $\alpha\hat{~}\la b'\ra$ will be eligible only finitely often.
With $n\in S$, this means that $P$ will be an infinite path through $T$,
picking out the least $b$ corresponding to each $a$ at the $\R_a$-node
$\alpha$, and picking out $\beta\hat{~}\la\cong\ra$ or
$\beta\hat{~}\la\fin\ra$ above a $\C_e$-node $\beta$
according as $C_e\cong F$ or not.

Now the list of fields $C_e$ includes every computable presentation
of every computably presentable field.  So, if $F$ is isomorphic to 
an arbitrary computable field $E$ (via an isomorphism $f$, say),
then that $E$ is precisely equal to some $C_e$.  We claim that
the $\C_e$-node $\beta$ on $P$ allows us to compute an isomorphism $g$
from $F$ onto $C_e$.  First, let $s_0$ be a stage after which $\beta$
is never initialized (so that no node to the left of $\beta$ is ever again
eligible).  Now for every $\R$-node $\alpha\subset\beta$, fix $b_\alpha\in\omega$
such that $\alpha\hat{~}\la b_\alpha\ra\subseteq\beta$.
Each of these $\alpha\hat{~}\la b_\alpha\ra$ is initialized
only finitely often, and the construction makes it clear that
each one, after its final initialization, adjoins only finitely many
elements to $F$: at most two witnesses and one tag.  Therefore,
there exists a stage $s_1\geq s_0$ after which no $\alpha\subset\beta$
ever again adjoins any elements to $F$.  Since the field $F_{s_1}$
is finitely generated, $f\res F_{s_1}$ is computable from the images
of its generators, which constitute finitely much information.
Hence we may set $g\res F_{s_1}=f\res F_{s_1}$.

It remains to define $g$ on elements adjoined by other $\R$-nodes $\alpha$.
If $\alpha$ lies to the right of $\beta$, then whenever $\alpha$ adjoins
any element to $F$ at some stage $s$ in its strategy for some $b$,
we simply wait until the next stage at which $\alpha$
is initialized.  Once this stage is complete, $\alpha$ never again adjoins
any elements from the splitting field of $p_{\alpha,b,s}(X)$,
and so once that stage is reached, we may find images for these
elements in $C_e$ (since $C_e\cong F$) and define $g$ to map them there.
(Of course, this uses Proposition \ref{prop:Ken}
and the choice of the $p$-polynomials to show
that every such splitting field is linearly disjoint from the compositum
of all the others, and that therefore these values for $g$ do not
interfere with the construction of $g$ on any other splitting field.)  

Finally, suppose $\beta\subset\alpha$.  Of course we do not know
whether such an $\alpha$ lies on $P$ or to its left or right.
However, when that $\alpha$ adjoins its first witness
$(x_{\alpha,b,s}+y_{\alpha,b,s})$ to $F$ at some stage $s$
for the $\alpha$-strategy for some $b$, we simply look for
the first root of $q_{\alpha,b,s}(X)$ to appear in $C_e$,
and let $g$ map the first witness to that root.  (Since $C_e\cong F$,
such a root must eventually appear in $C_e$, and by linear disjointness,
this extension of $g$ is still a field embedding.)
If the $\alpha$-strategy for $b$ never moves beyond
Step 1, then $F$ contains no more elements of the splitting field
of $p_{\alpha,b,s}(X)$, and so this is sufficient.
If it continues to Step 2 and adds the tag $x_{\alpha,b,s'}$
to $F$ at some stage $s'>s$, then we wait for such a tag
to appear in $C_e$ and define it to be $g(x_{\alpha,b,s'})$.
Notice that even if $\alpha$ eventually adjoins a second
witness to $F$ at a later stage $s''$, the first witness
to appear in $C_e$ \emph{must} be the one with the tag.
This follows from Step 3 of the construction for $\R$-nodes,
in which $\alpha$ waits until $C_e$ contains exactly one witness
node and also contains a tag for that node.  If $C_e$ acquired
a second witness before it acquired the tag for the first one,
then the construction would never have adjoined the second
witness to $F$, and $C_e$ would not be isomorphic to $F$,
contrary to hypothesis.  So $C_e$ must have produced the
tag for $g(x_{\alpha,b,s}+y_{\alpha,b,s})$ before
adjoining any second witness, and therefore it was safe
for us to define $g$ as we did on the first witness in $F$.
When (and if) $F$ acquires a tag for its first witness (in Step 2),
$C_e$ must subsequently acquire a tag for its
own first witness (in order to be isomorphic to $F$),
and then the second witness $(u_{\alpha,b,s''}+v_{\alpha,b,s''})$
to appear in $F$ (if $\alpha$ should execute Step 3 in its strategy)
will be matched by an (untagged) witness in $C_e$, to which
$g$ maps the second witness in $F$.  Thus we can compute
the value of this $g$ on every generator of $F$, and so $g$
is a computable field embedding of $F$ into $C_e$.  But
with $C_e\cong F$ by assumption, Corollary \ref{cor:endo}
shows that this $g$ is then an isomorphism.  Hence $F$
is computably categorical.

Next, suppose that $n\notin S$, and fix the least $a$ such that
no $b$ satisfies $f(n,a,b)\in\Inf$.  Now as argued above,
each node on the true path $P$ at any level $\leq 2a$
will have a successor on $P$.  When we reach the
$\R_a$-nodes at level $2a+1$, however, the $\alpha\in P$
at that level will have no successor eligible infinitely often,
since $(\forall b)f(n,a,b)\notin\Inf$.  We claim that instead,
the field $E_\alpha$ built by this $\alpha$ after its last initialization
is isomorphic to $F$, yet not computably isomorphic to $F$.
Since $E_\alpha$ is clearly a computable field (given finitely
much information, namely the last stage at which $\alpha$ was initialized),
this will show that $F$ is not computably categorical.

To see that $F\cong E_\alpha$, we begin at the first stage $s_0$
at which $\alpha$ is eligible after its last initialization.  At this stage
$E_{\alpha,s_0}$ is defined to be $F_{s_0}$ itself.
At all subsequent stages, the construction
(for every node $\alpha'$, not just $\alpha$) never
adjoins an element to $F$ without adjoining a corresponding
element to $E_\alpha$.  The only exceptions to this rule
are performed by $\alpha$ itself, at Step 2 of its strategies for
various values of $b$:  in Step 2 at those stages $s$,
$\alpha$ adjoins $x_{\alpha,b,s}$ to $F$ (which already
contained the witness $(x_{\alpha,b,s}+y_{\alpha,b,s})$)
without adjoining any element to $E_\alpha$ (which already
contained a witness element $(\xtilde_{\alpha,b,s}+\ytilde_{\alpha,b,s})$
of its own).  But at all subsequent stages, $\alpha$
will attempt to execute Step 3 for this $b$.  It will not
be allowed to do so as long as any $\C_e$-node
$\beta$ with $\beta\hat{~}\la\cong\ra\subseteq\alpha$
prevents it, which occurs if that $C_e$ fails to contain
exactly one witness for the $\alpha$-strategy for $b$,
along with a tag for that witness.  However, if this $C_e$
prevented it forever in this manner, then $C_e$ would
not be isomorphic to $F$, contradicting the fact that
such a $\beta\hat{~}\la\cong\ra$ must lie on $P$.
Therefore, eventually each of the finitely many $\C$-nodes
below $\alpha$ gives permission for $\alpha$ to execute
Step 3 in its strategy for $b$.  In doing so, $\alpha$
adjoins to $E_\alpha$ a new tagged witness, and adjoins
to $F$ a new untagged witness.  Moreover, by linear
disjointness, no more elements of the splitting field
of $p_{\alpha,b,s}(X)$ ever again enter either $F$
or $E_\alpha$,  Thus the witnesses and
tags in $E_\alpha$ and $F$ can be paired up perfectly,
and so indeed $E_\alpha$ is isomorphic to $F$.

Finally, suppose that some $\phi_b$ were an isomorphism
from $F$ onto $E_\alpha$.  Then, at some stage $t_b$
after which $W_{f(n,a,b)}$ receives no more elements,
the construction will have adjoined a first witness element
$(x_{\alpha,b,t_b}+y_{\alpha,b,t_b})$ to $F$ for $b$.
The isomorphism $\phi_b$ must map it to the corresponding witness
$(\xtilde_{\alpha,b,t_b}+\ytilde_{\alpha,b,t_b})$
adjoined to $E_\alpha$ at the same stage, since these
elements have no other conjugates in their fields at
that stage, and none are ever added unless $\phi_b$
maps the witness in $F$ to that in $E_\alpha$.  But once it does,
$\alpha$ executes Step 2, adjoining a tag for the witness in $F$,
and then (as we saw just above) eventually executes Step 3
and adjoins a new tagged witness in $E_\alpha$ and a new
untagged witness in $F$.  Therefore, $\phi_b$ maps
the tagged witness $(x_{\alpha,b,t_b}+y_{\alpha,b,t_b})$ in $F$
to the untagged witness
$(\xtilde_{\alpha,b,t_b}+\ytilde_{\alpha,b,t_b})$ in $E_\alpha$,
and so $\phi_b$ is not an isomorphism after all.  Since this holds
for every $b$, $F$ is not computably categorical.  This completes the
proof of Theorem \ref{thm:Pi04complete}.
\qed\end{proof}

At first glance, the foregoing proof appears to be a standard $\emptyset''$
construction, using the true path $P$ through a computable tree.
However, a $\emptyset''$ oracle is not in fact enough to compute $P$.
It can compute the successor on $P$ of any $\C$-node $\beta\in P$,
and it can compute the successor of an $\R$-node $\alpha\in P$
provided that $\alpha$ has one.  However, $P$ may actually end
at $\alpha$ (in which case $E_\alpha$ is the computable field
showing that $F$ is not computably categorical), and this situation holds
if and only if $\forall b f(n,a,b)\notin\Inf$, which is a $\Pi^0_3$ condition.
So in fact, to compute $P$ and recognize when it terminates (if ever),
a $\emptyset'''$ oracle is required.

\section{Conclusions and Questions}
\label{sec:conc}

The ultimate goal of this project was to provide a structural
characterization for computable categoricity for algebraic fields.
The main question, therefore, is the extent to which we have achieved
this goal.  Admittedly, the goal itself is somewhat vague:
what constitutes a structural characterization?  A first-order
property in model theory would be the ideal result, but this goal
seems beyond reach.

For illumination on this question, consider the characterization
of computably categorical linear orders $\L$ as those with only
finitely many adjacencies.
This property is not expressible in first-order languages,
as one quickly proves using the Compactness Theorem.
It is also readily seen to be a $\Sigma^0_3$-complete
property, and so, in terms of complexity, we know exactly
the level of difficulty of deciding computable categoricity
for computable linear orders.

Notice also that, because computable categoricity implies
relative computable categoricity for linear orders, another
equivalent characterization would be the existence of a
$\Sigma^0_1$ Scott family for $\L$.  This property
is also $\Sigma^0_3$-complete, for linear orders
as for computable structures in general, and could
also be taken as a characterization of computable categoricity.
However, it is vastly less satisfying than the characterization
by the number of adjacencies:  the latter feels much more
``structural.''  To quantify this, we note that the characterization
using adjacencies can be expressed as a computable
$L_{\omega_1\omega}$ formula (that is, with countable
conjunctions and disjunctions allowed) in the language
of linear orders.  In addition, the proof of the equivalence
of the latter to computable categoricity makes it clear exactly
how the property of finitely many adjacencies corresponds
to computable categoricity, much more clear than can be said
of the characterization by Scott families.  So we consider the
characterization by adjacencies to be the better characterization.

Since the initial consideration of computable categoricity for fields
by Fr\"ohlich and Shepherdson in \cite{FS56}, the problem of
characterizing computable categoricity for fields has not given much ground.
Without offering specific justification,
we suspect that the results in this article are as good as one is likely
to get in the case of algebraic fields.  As far as complexity, that
statement can be quantified:
$\Pi^0_4$-completeness of computable categoricity for algebraic fields,
demonstrated in Theorem \ref{thm:Pi04complete}, pinpoints the complexity
of the notion.  Likewise, of course, the characterization for relative computable
categoricity turned out to be $\Sigma^0_3$ (as it must, by the
work in \cite{AKMS89} and \cite{Chis})
and complete at that level (as commonly happens for relative
computable categoricity).  As usual, the characterization by Scott families
is unsatisfying, and we consider Theorem \ref{thm:rcc} to be a significant
step forward, since it equates this characterization to the more structural
notion in items (4) and (5) of that theorem.
It is not clear that any more satisfactory characterization
of relative computable categoricity is likely to be discovered.

For computable categoricity, we likewise consider Theorem
\ref{thm:Pi04complete} to be substantial progress.  Nevertheless,
the result still feels less satisfactory.  The property given in Proposition
\ref{prop:pi4} is really just the definition of computable categoricity,
in the specific context of algebraic fields.  Theorem \ref{thm:Pi04complete}
then shows that one cannot do better, in terms of complexity,
and we consider it important to recognize that in this context,
Definition \ref{defn:cc} can achieve the minimum possible complexity,
simply by replacing the notion of classical isomorphism
by an equivalent statement (namely the condition from
Corollary \ref{cor:iso}).  We believe that this is the
first known instance of this phenomenon.
However, it still does not seem impossible that a ``more structural''
characterization might be found.

We attach additional importance to Theorem \ref{thm:Pi04complete} because
of the new level of complexity it exhibits.  Previous characterizations
of computable categoricity for standard classes of computable structures
have generally shown it to be $\Sigma^0_3$-complete
(and equivalent to relative computable categoricity):
this situation holds for linear orders, Boolean algebras, trees (as partial
orders), and ordered abelian groups, for example.
Relative computable categoricity is widely viewed as a ``nicer''
property, largely because of its straightforward syntactic characterization
in \cite{AKMS89} and \cite{Chis}, and it was already known that
computable categoricity has strictly higher complexity than relative
computable categoricity in many well-known classes of structures,
such as graphs, partial orders, groups, and rings.
In \cite{W03}, White showed that for computable graphs,
computable categoricity is $\Pi^0_4$-hard, and
\cite{HKSS02} allows the complexity result to be carried over to
the other well-known classes mentioned (although it only proves
computable categoricity to be $\Pi^0_4$-hard in those classes,
not necessarily $\Pi^0_4$-complete).
The fact that computable categoricity turned out to be $\Pi^0_4$-complete
for algebraic fields took us rather by surprise, as this is the first everyday
class of mathematical structures in which it turned out to be a
$\Pi^0_n$-complete property (as opposed to $\Sigma^0_n$-complete)
for any $n$ at all.  Indeed, to our knowledge, algebraic fields are
the first standard class of structures for which the complexity of
computable categoricity has been determined and has turned out
not to be $\Sigma^0_3$-complete.

(For careful readers, we point out a small error in the
final paragraph of \cite{W03}, where it is asserted that computable
categoricity is $\Pi^0_3$-complete for the class of algebraically closed fields.
In fact, for such fields, Ershov \cite{E77} showed it to be equivalent
to the property of having finite transcendence degree,
which is $\Sigma^0_3$-complete and is also equivalent
to relative computable categoricity for such fields.
Likewise, as of the writing of \cite{W03}, all other known
index sets for computable categoricity were $\Sigma^0_3$,
not $\Pi^0_3$ as stated there.)

It should be noted that the class of algebraic fields is not
first-order definable:  every axiom set that holds in all
algebraic fields will hold in certain non-algebraic fields as well.
This fact might help explain the unusual level of complexity.
In characteristic $0$, our theorems carry over to fields
of finite transcendence degree over $\Q$, since essentially
all constructions can be carried out after replacing $\Q$
by a purely transcendental subfield $\Q(X_1,\ldots,X_k)$
over which $F$ is algebraic.  (Alternatively, for transcendence
degree $k$, just enrich the signature by $k$ constants,
with axioms saying that they are algebraically independent over $\Q$.)
For fields of infinite transcendence degree, the question of
computable categoricity is not trivial: most such fields are not
computably categorical, but the work of Miller and Schoutens
in \cite{MS08} proved the existence of a computably
categorical field of infinite transcendence degree.  One would guess
that computable categoricity has even higher complexity
for the class of all fields; it certainly cannot become any lower
than $\Pi^0_4$, since algebraic fields form a subclass.

Finally, we mention computable dimension for algebraic fields.
Goncharov defined the \emph{computable dimension} of a computable
structure to be the number of computable presentations of that structure,
up to computable isomorphism.  He showed that every cardinal
from $1$ through $\omega$ can be the computable dimension of a computable
structure.  (See \cite{G82} and \cite{G98} for these and related results.)
However, by far the most common computable dimensions
are $1$ (which is equivalent to computable categoricity) and $\omega$,
and for many classes of structures, these are the only
possible computable dimensions:  linear orders, Boolean algebras,
and trees, for example.  At one time, we believed that we had proven
this to be true of algebraic fields as well.  However, a problem subsequently
was found in the proof of the theorem (from another article) which we used
to show the impossibility of finite computable dimension $>1$.
A recent result in \cite{MPPSS14} shows that in computable fields
of infinite transcendence degree over $\Q$, all computable dimensions
are possible.
However, the question of finite computable dimension for computable
algebraic fields remains open.  We conjecture it to be impossible
in purely quadratic algebraic extensions of $\Q$ (i.e., those algebraic
extensions $F$ such that for every finitely generated subfield $E$ of $F$,
the degree of $E$ over $\Q$ is a power of $2$), and also in the related class
containing all computable finite-branching subtrees of $\omega^{<\omega}$.
It is known that the latter class includes structures which are computably
categorical but not relatively so (by a proof very similar to that in
Section \ref{sec:cc}), and so this conjecture would imply that the
condition of computable categoricity without relative computable categoricity
need not entail the existence of structures of finite computable dimension $>1$:
the one pathology can occur without the other.

\bibliographystyle{amsplain}

\end{document}